\documentclass{amsart}

\usepackage{amsmath, amsthm, amssymb}
\usepackage{graphicx}

\newtheorem{theorem}{Theorem}[section]
\newtheorem{lemma}[theorem]{Lemma}
\newtheorem{Prop}[theorem]{Proposition}
\newtheorem{Cor}[theorem]{Corollary}

\theoremstyle{definition}

\newtheorem{example}[theorem]{Example}

\theoremstyle{remark}
\newtheorem{remark}[theorem]{Remark}
\numberwithin{equation}{section}

\newcommand{\End}{\operatorname{End}}
\newcommand{\modc}{\operatorname{-mod}}

\newcommand{\Span}{\operatorname{span}}

\newcommand{\Res}{\operatorname{res}}
\newcommand{\Mod}{\operatorname{mod}}

\newcommand{\R}{{\mathbb R}}
\newcommand{\C}{{\mathbb C}}
\newcommand{\Z}{{\mathbb Z}}

\newcommand{\A}{{\mathcal A}}
\newcommand{\Bip}{\operatorname{Bip}}
\newcommand{\KBip}{\operatorname{KBip}}

\newcommand{\Par}{\operatorname{Par_{\emptyset}}}
\newcommand{\rank}{\operatorname{rank}}
\newcommand{\Q}{{\mathbf Q}}
\newcommand{\q}{{\mathbf q}}

\begin{document}    

\title{Cell structures on the blob algebra } 
\author[S. Ryom-Hansen] {Steen Ryom-Hansen}
\address{Instituto de Matem\'atica y F\'isica, Universidad de Talca \\Chile }
\email{steen@inst-mat.utalca.cl }
\thanks{Supported in part by FONDECYT grants 109070 and 1121129, by Programa Reticulados y Simetr\'ia
and by the MathAmSud project OPECSHA 01-math-10.}

\subjclass[2010]{20G05, 20C08, 05E10}

\begin{abstract}
{We consider the $ r = 0 $ case of the conjectures
by Bonnaf\'e, Geck, Iancu and Lam on 
cellular structures on the Hecke algebra of type 
$ B $. We show that this case induces the natural cell structure on the 
blob algebra $ b_n $ by restriction to one-line bipartitions.}
\end{abstract}

\maketitle

\section{Introduction}
The purpose of this article is to continue the investigation, 
initiated in [RH], 
of the relationship between the 
representation theories of the Hecke algebra 
${\mathcal H}_n = {\mathcal H}_n(Q,q) $ of type $ B $ 
and of the blob algebra $ b_n = b_n(q,m) $. 
The Hecke algebra ${\mathcal H}_n $ of type $ B $ is a well-known two-parameter 
deformation of the hyperoctahedral group
whereas 
the blob algebra $b_n $, introduced in [MS] from motivations in statistical mechanics, is a 
diagram algebra of marked (blobbed) Temperley-Lieb 
diagrams. A main point of our work, already present in [RH], is that $ b_n $ can also be realized as a 
quotient of $ {\mathcal H}_n $ thus making the $ b_n $-representations $ {\mathcal H}_n $-representations 
by inflation. Viewing $ b_n $ as a quotient of $ {\mathcal H}_n $ is analogous to viewing 
the Temperley-Lieb algebra $ TL_n $ as a quotient of the Hecke algebra of type $ A $, and indeed 
$ b_n $ is also sometimes called the Temperley-Lieb algebra of type $ B$.  

\medskip
Dipper-James-Murphy introduced in [DJM] for each bipartition $ (\lambda , \mu ) $ of total degree $ n $ a Specht 
module $ S_n(\lambda, \mu) $ for $ {\mathcal H}_n $.
Let $ J_n $ be the kernel of the quotient map $ {\mathcal H}_n  \rightarrow b_n $. 
We then showed in [RH] that $ J_n  S_n(\lambda, \mu) = 0 $ as long as $ (\lambda, \mu )$ is a one-line 
bipartition 
and so these $ S_n(\lambda, \mu) $ factor over the quotient map to become $ b_n $-modules. 
One might now suspect that $ S_n(\lambda, \mu )$ is a standard module for the quasi-hereditary 
algebra $ b_n $. Indeed, we showed 
that many properties of the standard modules
are shared by the $ S_n(\lambda, \mu )$, but somewhat surprisingly we could prove in [RH]
that they do not verify the 
relevant universal property and so do not identify with standard 
modules, except in trivial cases. 

\medskip
Recall G. Lusztig's monograph [Lu2] on the representation theory of Hecke algebras 
with unequal parameters.
Let $ W $ be a Coxeter group and let $ L: W \rightarrow \Gamma $ be a length function in the 
sense of [Lu2], for 
$ \Gamma $
a totally ordered Abelian group. Associated to this data, [Lu2]
contains a construction of cells in $ W$ and cell modules for the corresponding multiparameter
Hecke algebra, generalizing the 
construction
from [KL] in the one-parameter case. When $ W $ is of type $ B$ 
the length function is specified by 
$ a:=\log q, b:=\log Q \in \Gamma $.
In [BGIL] a series a conjectures were formulated for type $ B $ which, if true, would 
put a high degree of structure on this. 
Assume that $ b \not\in \{ a, 2a, \ldots , (n-1) a \} $ 
and that $ a $ and $ b $ are positive in $ \Gamma$.
According to the conjectures,
the setting should 
give rise to a {\it cellular algebra} datum on $ {\mathcal H}_n $ in the sense of Graham and Lehrer, 
where the underlying poset $ \Lambda $ should be the set of bipartitions $ \Bip(n) $ of total degree $ n $ 
with partial order and  
map $ \Lambda \times \Lambda \rightarrow {\mathcal H}_n  $ defined in terms of a certain 
{\it domino insertion} algorithm, depending on $ \Gamma $. Furthermore, 
by the work of Bonnaf\'e and Jacon [BJ], 
the different cellular algebra structures on $ {\mathcal H}_n $ should account for
the different ways of parameterizing the simple modules for $ {\mathcal H}_n $ that are given by Ariki's Theorem in [A].

\medskip
These conjectures have only been fully proved in the so-called asymptotic case $ b > (n-1)a $, 
see [BI], where the cell modules turn out to be the ones given by Dipper-James-Murphy.
In this work we focus on the case $ \Gamma := \Z $, $ a:= 2 $ and $ b = 1 $. 
This is another extreme case since 
$ b < a $ and so $ r=0 $ in the [BGIL] notation. 
We show that the poset structure on $ \Bip(n) $ in this case is compatible with the quasi-hereditary order on the
category of $ b_n$-modules when restricted to one-line bipartitions, 
the map being given by $ (\lambda, \mu )   \mapsto k-l  $ where $ \lambda = (k) $ and 
$ \mu = (l) $. We show that the ideal $ J_n $ is generated by the set of 
Kazhdan-Lusztig elements $ C_w $ for which $ w $ does not correspond to a one-line bipartition. We moreover 
show that the cell module given by the one-line bipartition $ (\lambda, \mu ) $
is isomorphic 
to the $ b_n $ standard module $ \Delta_n(k-l) $
where $ \lambda = (k) $ and 
$ \mu = (l) $. 
To summarize our findings: the $ a= 2, b = 1  $ 
case of the [BGIL] conjectures induces the blob algebra category when 
restricted to one-line bipartitions.

\medskip
This given, the algorithm described in [Ja] can be used to answer the question that was raised in 
[RH], namely to describe the
Kleshchev bipartition that corresponds to the simple $ b_n $-module 
$ L_n(\lambda) $. 

\medskip
Let us indicate the layout of the article. 
The first section contains a combinatorial analysis of the domino insertion 
algorithm mentioned above. The main result is a characterization
of the elements $ W_b$ of 
the Weyl group $ W_n $ of type $ B $ that go to two-line tableaux under domino insertion.
This characterization uses the Coxeter presentation of $ W_n $. The section relies on results
of Taskin, [T].

\medskip
In the next section we recall the presentation of $ b_n $ 
as a quotient of $ {\mathcal H}_n $ and show that the defining ideal is given by the Kazhdan-Lusztig type 
elements $ C_w \in {\mathcal H}_n $ where $ w \notin W_b$.
In the following section  
we show our main results, identifying the cell modules with the standard modules. 
To be more precise, we show that the cell modules verify the 
universal property for the standard modules, given within the 
framework of the globalization-localization formalism.
For this to work we rely on Lusztig's results in [Lu1] that we combine 
with the results of Fan and Green [FG] on type $ A $.

\medskip
Finally, in the last section we show how the Fock space approach to the representation theory of 
$ {\mathcal H}_n $ can be used to reprove the main results of [MW] and to obtain the 
Kleshchev bipartitions of the simple modules for $ b_n $.

\medskip
It is a great pleasure to thank the referees for many useful comments and suggestions.

\section{Basic notation and domino insertion}
In this section we first fix  
some notation that shall be used throughout the article. 
We then investigate the 
domino insertion algorithm 
for the Weyl group of type $ B $. We describe the elements that are mapped under it to two-line partitions, that is 
domino tableaux whose underlying partition has at most two lines (parts). 

\medskip
We shall throughout assume knowledge of the definition and 
basic properties of the Robinson-Schensted algorithm, as
exposed in for example [Sa]. 

\medskip
For the following basic combinatorial concepts related to the Weyl group of type $ B $, we refer the reader to
section 8.1 of [BB]. 
Let $ W_n $ be the Weyl group of type $ B_{n} $. 
It is a Coxeter group on generators 
$ s_0,  s_1, \ldots, s_{n-1} $ with relations
$$ \begin{array}{lr}   s_i^2 = 1 &  \mbox{  for } i = 0, \ldots , n-1 \\ 
(s_i  s_{i+1})^3 = 1 & \mbox{for } i = 1, \ldots , n-2   \\
(s_i  s_j)^2 = 1      & \mbox{ for } | i-j | > 2   \\ (s_0  s_1 )^4 = 1.  & \\
\end{array} 
$$
Let $ I_n:= I_n^+ \cup I_n^- $ where $ I_n^+ := \{  1,  2, \ldots,  n \} $ and 
$ I_n^- := \{  -1,  -2,  \ldots,  -n \} $. Then $ W_n $ can also 
be described as the subgroup of the symmetric group on the 
elements $ I_n$ 
generated by $ s_0 := (-1,1) $ and 
$$ s_i := (i,i+1)(-i,-i-1) $$ in cycle notation. We shall adopt the convention that cycles 
are multiplied from right to left. 
The subgroup of $ W_n $ generated by 
$ s_1, s_2, \ldots ,s_{n-1} $ is the symmetric group $ S_n $.

\medskip
For elements $ w \in W_n $ we shall also use {\it word} or {\it sequence} notation as follows
$$ w = i_1 i_2 i_3 \ldots i_n $$
where $ i_k \in I_n $. By this we mean that $ w  $ acts (on the left) on $ I_n $ as follows
$$ w: 1 \mapsto i_1, 2 \mapsto i_2, \, \ldots,  n \mapsto i_n $$ 
and then also necessarily $ -1 \mapsto -i_1,  -2 \mapsto -i_2, \ldots,  -n \mapsto -i_n$.
In this setting we use the standard notation $ \overline{i} := -i \in I_n^- $ for $ i \in I_n^+$. 
Thus $ i $ appears in $w=i_1 i_2 i_3 \ldots i_n \in W_n $ if and only if $ \overline{i} $ does not not appear.	

\medskip
It is normally clear whether 
a given $ w \in W_n $ is written as a product of Coxeter generators
or as a word over $I_n $ and we shall therefore generally not explicitly mention 
the chosen form.
If for example $ w := s_0  s_1  s_2 \in W_3 $, we may write 
$$  w = s_0  s_1  s_2= 2 \, 3 \, \overline{1}.   $$

\medskip
We denote by $ < $ the Bruhat-Chevalley 
order on $ W_n$ where by convention the neutral element $ 1 \in W_n $ is the smallest of all.
Assume that  $ w = i_1 \, i_2 \, i_3 \ldots \, i_n  \in W_n $. Then the 
following conditions describe the right descent set of $ w $ with respect to $ < $
$$ \begin{array}{ll} w  s_k < w  \mbox{ iff } i_k > i_{k+1} & \mbox{ for } k = 1, 2,  \ldots , n-1 \\  
w  s_0 < w \mbox{ iff } i_1 < 0. & \end{array} $$
If $ w \in W_n $ is written in word form, its right descent set can be used to write it 
as a reduced expression in the Coxeter generators $ s_i $. 
\begin{example} Assume that $ w = 3 \, \overline{1} \, \overline{2} \, 4 $. Then $ s_1  s_0  s_1 
s_0  s_2  s_1  $ is a reduced expression for $ w $ obtained from the above description of 
the right descent set. Indeed, $ s_2  s_1 $ moves $ 3 $ past $ \overline{1} \, \overline{2} $, then 
$ s_0 $ changes $ \overline{1} $ to $ 1 $ and finally $ s_1  s_0  s_1 $ changes $ \overline{2} $ to 
$ 2 $. 
\end{example}
Throughout the paper, we shall be specially interested in the subset $ W_b = W_{n,b} $ of $ W_n $. It consists of 
those $ w \in W_n $ that have no reduced expressions
$ w = s_{i_1}  s_{i_2}  s_{i_3} \ldots  s_{i_N} $
that contain a subexpression $ s_{i_k}  s_{i_{k+1}}  s_{i_{k+2}} $ of the form 
$$ \begin{array}{lr} s_i  s_{i \pm 1 }  s_i \, \,  \, \mbox{ for } i = 1, 2, \, \ldots , n-2 & 
\mbox{ or     } s_{n-1}  s_{n-2}  s_{n-1}.
\end{array}
$$
Thus the subexpression $ s_0  s_1  s_0 $ is allowed whereas $ s_1  s_0  s_1 $
is not. 

\medskip
Our aim is to describe the image of $ W_b $ under the domino insertion correspondence described 
for example in [BGIL]. In order to do so we first need a description of $ W_b $ in terms of words.
This description will only be indirect, but for our purposes this will be sufficient. 

\medskip
In general we use the convention that empty index sequences correspond to the void subsequence. 
For example, in the next Lemma, the case $ k = 1 $ corresponds to $ w= \overline{a}_1 \,  i_{2 } \ldots i_n $.

\begin{lemma}{\label{first_lemma}} 
Assume that $ w \in W_n $ and assume that it can be written as follows
$$ w = i_1 \, i_2 \,  \, \ldots \, i_{k-1} \, \overline{a}_1 \,  i_{k+1 } \ldots i_n $$
where $ i_1, i_2 , \ldots i_{k-1}, a_1 > 0 $.
Then $ w \in W_b $ if and only if
$$ \begin{array}{lr} a_1 < i_1 < i_2 < \ldots < i_{k-1} & \mbox{ and   } \, \, 
 a_1 \,i_1 \, i_2 \, \ldots \, i_{k-1}  \,  i_{k+1 } \ldots i_n \in W_b. \end{array} $$
\end{lemma}
\begin{proof}
Suppose that $ w \in W_b $ and define 
$ w_1 \in W_n $ as
$$ w_1 = a_1 \, 
i_1 \, i_2  \, \ldots \, i_{k-1}  \, i_{k+1} \,  \ldots i_n. $$
Using the above description of the right descent set we get that $ w $ has a reduced expression of the form
\begin{equation}\tag{$\ast$} w = w_1 s_0  s_{1} \ldots  s_{k-2} s_{k-1}    \end{equation}
and the second statement follows, since any reduced expression for $ w_1 $ can be extended 
to a reduced expression for $w$. 

\medskip
If now $ a_1 < i_1 < i_2 < \ldots < i_{k-1} $ is not satisfied then by the 
description of the right descent set there will be an index 
$1 \leq  j \leq k-1 $ such that $ w_1 s_j < w_1 $. But by 
formula ($\ast$)
this contradicts the assumption that $ w \in W_b $.

\medskip 
To show the other implication we assume 
that $ a_1 < i_1  < \ldots < i_{k-1} $ holds, that 
$ w_1 = a_1 \, 
i_1 \, i_2  \, \ldots \, i_{k-1}  \, i_{k+1} \,  \ldots i_n \in W_b $ and 
that $ w \not\in W_b$.
Since $ s_0 s_1 s_2 \ldots s_{k-1} $ is a  unique representation for $ w_1^{-1} w $ and since 
$ w_1 \in W_b $ we conclude that $ w_1 $ must have a reduced expression of the form $ w_1 := w_2 s_j $, 
for an index $ j $ such that $ 0 \leq j \leq k-1 $. But then $ s_j $ belongs to the right descent set for
$ w_1 $, contradiction.
\end{proof}

Before stating our next result, we need to recall the combinatorial notion of a
{\it decreasing subsequence}.
Suppose that $ w = i_1 i_2 i_3 \ldots i_n \in W_n $. 
A decreasing subsequence of $ w $ of length $ k $ is defined to be 
a subsequence $ i_{\iota_1} i_{\iota_2} \ldots i_{\iota_k} $ of $ w $ such that $ \iota_j < \iota_{j+1} $
and $ i_{\iota_j} > i_{\iota_{j+1}} $ for $ j = 1, \ldots, k $.

\medskip
Setting $W_c:=W_b \cap S_n $, it is 
known that $W_c$ can be described as the words over $ I_n^+ $
with no decreasing subsequences of length strictly greater than two
and so it 
corresponds under the Robinson-Schensted algorithm  to pairs of 
tableaux $ (s, t ) $ of a two-line partition $ \lambda=(\lambda_1, \lambda_2)$.

\begin{theorem}{\label{first_thm}} 
Suppose $ w \in W_n $ and write it as 
$$ w = i_1 \, i_2 \, \, \ldots \, i_{k_1-1} \, \overline{a}_1 \, i_{k_1+1} \,   \ldots  \, i_{k_2-1}
\,  \overline{a}_2 \, 
i_{k_2+1} \, 
\ldots i_{k_l-1} \, \overline{a}_l \,  
i_{k_l+1} \, \ldots \,  
i_n $$
where $ \overline{a}_1, \ldots \overline{a}_l $ are the only negative numbers that occur in $ w$.
Define $$ w^l := a_l \, a_{l-1} \, \ldots \, a_1 \, i_1 \, i_2   \ldots i_n. $$
Then $ w \in W_b $ if and only if 
$$a_l < a_{l-1 }  < \ldots < a_1 < i_1 < i_2  < \ldots < i_{k_l -1} $$ and 
$ w^l $ has no decreasing subsequences of length strictly greater than $2$. 
\end{theorem}
\begin{proof}
Suppose first that $ w \in W_b $. We generalize $ w^l $ as follows 
$$ \begin{array}{l}
w^1 = a_1
i_1 \, i_2 \, \, \ldots \, i_{k_1-1}  \, i_{k_1+1} \,   \ldots  \, i_{k_2-1}
\,  \overline{a}_2 \, 
i_{k_2+1} \, 
\ldots i_{k_l-1} \, \overline{a}_l \,  
i_{k_l+1} \, \ldots \,  
i_n  \\
w^2 = a_2 \, a_1
i_1 \, i_2 \, \, \ldots \, i_{k_1-1}  \, i_{k_1+1} \,   \ldots  \, i_{k_2-1}
\,  
i_{k_2+1} \, 
\ldots i_{k_l -1} \, \overline{a}_l \,  
i_{k_l+1} \, \ldots \,  
i_n \\ \vdots \\
w^l = a_l a_{l-1} \ldots  a_1
i_1 \, i_2 \, \, \ldots \, i_{k_1-1}  \, i_{k_1+1} \,   \ldots  \, i_{k_2-1}
\,  
i_{k_2+1} \, \ldots \, 
 i_{k_l-1} \, 
i_{k_l+1} \, \ldots \,  
i_n.
\end{array}
$$
By the proof of the previous Lemma we have $ w^k \in W_b $ for all $ k $ and so we get the inequalities 
$$a_l < a_{l-1 }  < \ldots < a_1 < i_1 < i_2  < \ldots < i_{k_l-1} $$ 
by using the previous Lemma recursively. 
But $ w^l \in W_b \cap S_n $ and so we have proved one implication of the Theorem.

\medskip
The other implication follows in a similar way from the previous Lemma.
\end{proof}

\begin{example} 
Let us consider $ w = 3 \, \overline{1} \, \overline{2} \, 4 $ from the previous example. 
Then, in the notation of the Theorem, we have $ w^2 := 2 1 3 4 $ and so $ w \not\in W_b$, since $ 213 $ is 
not increasing.
The conclusion could also have been obtained directly from the definition of $ W_b $ and the description of $ w $ 
found in the previous example:
$ w = s_1  s_0  s_1 s_0  s_2  s_1  $. 
\end{example}

\begin{example} 
Using the Theorem, one can produce elements of $ W_b$ by shuffling an increasing
sequence of negative numbers with an increasing sequence of positive numbers, such that the positive  
terms all have absolute values larger than the negative numbers. For example
$$ w:=  \overline{4} \, \overline{3}\, 5 \, \overline{2} \, 6 \, 7 \, 8 \, \overline{1} \, 9  \in W_{9, b}.$$
\end{example}

\begin{figure}
\includegraphics{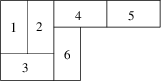}
\caption{}
\label{one}
\end{figure}

The notion of domino tableaux shall be important to us. A domino tableau is the Young diagram of an integer 
partition of $ 2n $ with node set partioned into dominoes, that is horizontally or vertically neighboring nodes. 
The dominoes are labeled with numbers $ 1, 2, \ldots, n $. A domino tableau is called standard if the labeling is 
increasing from left to the right and from top to bottom.
Let $ SDT(n) $ denote the set of standard domino tableaux in $ n $ dominoes. Figure \ref{one} gives 
an example from $ SDT(6)$.
We define $ SDT := \bigcup_n  SDT(n) $. For $ S \in STD $ we let $ Sh(S) $ denote the shape of 
its underlying partition. Let 
$ SDT^2(n) $ be the set 
$$ SDT^2(n) := \{\,  (S, T)  \in  SDT(n) \times SDT(n) \, | \,  Sh(S)= Sh( T)   \}. $$
The domino insertion algorithm
establishes a bijection between $ W_n $ and $ SDT^2(n) $.
It was introduced in [BV] as a generalization of the Robinson-Schensted algorithm to type $ B$.
A slightly different version of the algorithm, using a bumping procedure, was introduced in [Ga], see also [vL].
We shall not here give a precise description of the algorithm, but 
refer the reader to 
for instance [BGIL] or [La].

\medskip
Let us denote by $ (P(w), Q(w) ) $ the pair of domino tableaux associated with $ w \in W_n $
under domino insertion.
We say that $ w $ and $ w_1 $ belong to the same Knuth (plactic) class, or 
$ w \stackrel{p}{\sim} w_1 $,
if $ P(w) = P(w_1) $.  Dually, 
we say that $ w $ and $ w_1 $ belong to the same dual Knuth (coplactic) class, or 
$ w \stackrel{p^{\ast}}{\sim} w_1 $,
if $ Q(w) = Q(w_1) $.  
 
\medskip
Taskin considers in Definition 3.1 of [T] a set of generalizations of the Knuth relations, 
$ D_i^r, \, i = 1, 2, \ldots 5$, 
that generate the (co)plactic relations. We now explain these relations in the case 
that we need, $ r= 0 $, where they simplify somewhat.	 
The elements of $ W_n $ are always assumed to be in word form. 

\medskip	
Let $ f: I_n \rightarrow I_n $ be any bijection such that $ f(1) < f(2) < f(3) $. Then 
$ D_1^0 $ can be reformulated as the combination of the following two relations
\begin{equation}{\label{first_Knuth}}
   \cdots f(2)f(3)f(1) \cdots \,  \stackrel{KT}{\sim} \,  \cdots f(2)f(1)f(3) \cdots 
\end{equation}
\begin{equation}{\label{second_Knuth}}
 \cdots f(1)f(3)f(2) \cdots \, \stackrel{KT}{\sim} \, \cdots f(3)f(1)f(2) \cdots 
\end{equation}
where we use the convention that there are no changes of dotted elements.
The relation $ D_2^0 $ is void whereas the relation $ D_3^0 $ is the following one
\begin{equation}{\label{third_Knuth}}
    i_1 \, i_2 \cdots \,  \stackrel{KT}{\sim} \,   \overline{i_1} \, i_2  \cdots \, \, \, \,
\mbox{ if } |i_1 | > | i_2 |
\end{equation}
under a further condition on the dotted elements that we do not need to detail.

\medskip
The remaining two relations $ D_4^0 $ and $ D_5^0 $ are more complicated to express than 
the first ones. But since we are only considering the $ r= 0 $ case of Taskin's results, 
we may use a somewhat simplified notation.

Let us first consider $ D_4^0 $. 
Assume that $ \alpha $ and $ \alpha^{\prime} $ are elements of $ W_n $ that can be expressed as 
$ \alpha = u \ldots  $ and $ \alpha^{\prime}  = u^{\prime}  \ldots $ 
where 
$$ 
\begin{array}{l}
 u =   a_{1,1} b_{1,1} a_{2,2} a_{2,1}  b_{2,2} b_{2,1}  \ldots 
a_{k,k} a_{k,k-1} \ldots a_{k,1}  (b_{k,k} b_{k,k-1} \ldots b_{k,1})
a_{k+1,k}  \ldots a_{k+1,1} z \\
u^{\prime} =   a_{1,1} b_{1,1} a_{2,2} a_{2,1} b_{2,2} b_{2,1} \ldots 
(-b_{k,k})
a_{k,k} a_{k,k-1} \ldots a_{k,1} (b_{k,k-1} \ldots b_{k,1})
a_{k+1,k}  \ldots a_{k+1,1} z \\
\end{array}
$$
for some $ z \in I_n $ and $ k \geq 1$.
(Notice that there is no $ a_{k+1,k+1}$).
Suppose moreover that $\{ a_{i, j } \}_{i,j \geq 1}  $ and $\{ b_{i, j } \}_{i,j \geq 1}  $ 
satisfy 
$$ \begin{array}{l}
a_{i, j } > 0  \mbox{ and } b_{i, j } < 0  \, \mbox{ (or vice versa) for all} \,\,   i, j \geq 1 \\
| a_{i, j-1 } | < | a_{i, j } | < | a_{i+1, j } | \mbox{ and }
| b_{i, j-1 } | < | b_{i, j } | < | b_{i+1, j } | \\
| b_{i, i} | < | a_{i+1, i+1 } | < | b_{i+1, i+1 } | \mbox{ for all } i = 1, \ldots ,  k-1. 
\end{array}
$$
 Then 
$ D_4^0 $ is the relation 
\begin{equation}{\label{fourth_Knuth}}
   \alpha \stackrel{KT}{\sim} \alpha^{\prime} 
\end{equation}
under certain further conditions on $ z $
that we do not detail. 

\medskip
\noindent
Let us finally consider the relation $ D_5^0 $. Assume that $ \alpha $ and $ \alpha^{\prime} $ 
are elements of $ W_n $ that can be expressed as 
$ \alpha = u \ldots  $ and $ \alpha^{\prime}  = u^{\prime}  \ldots $ 
where 
$$ 
\begin{array}{l}
 u =   a_{1,1} b_{1,1}   \ldots 
(a_{k,k}  \ldots a_{k,1} ) (b_{k,k} \ldots b_{k,1}) 
(a_{k+1,k+1} a_{k+1,k} \ldots  a_{k+1,1} ) (b_{k+1,k} \ldots b_{k,1})
z \\
u^{\prime} =  
a_{1,1} b_{1,1}   \ldots 
(a_{k,k}  \ldots a_{k,1} ) (-a_{k+1,k+1})(b_{k,k} \ldots b_{k,1}) 
(a_{k+1,k}  \ldots  a_{k+1,1} ) (b_{k+1,k} \ldots b_{k,1}) z
\end{array}
$$
for $ z \in I_n$ and $ k \geq 1 $. (This time there is no $ b_{k+1,k+1}$).
Assume moreover that $\{ a_{i, j } \}_{i,j \geq 1}  $ and $\{ b_{i, j } \}_{i,j \geq 1}  $ 
satisfy 
$$ \begin{array}{l}
a_{i, j } > 0  \mbox{ and } b_{i, j } < 0  \, \mbox{ (or vice versa) for all} \,\,   i, j \geq 1 \\
| a_{i, j-1 } | < | a_{i, j } | < | a_{i+1, j } | \mbox{ and }
| b_{i, j-1 } | < | b_{i, j } | < | b_{i+1, j } | \\
| a_{i, i} | < | b_{i, i } | < | a_{i+1, i+1 } | \mbox{ for all } i = 1, \ldots, k. 
\end{array}
$$
Then 
$ D_5^0 $ is the relation that 
\begin{equation}{\label{fifth_Knuth}}
   \alpha \stackrel{KT}{\sim} \alpha^{\prime} 
\end{equation}
under certain further conditions on $ z $
that, once again, we do not detail. 

\medskip
We shall refer to the relations 
(\ref{first_Knuth}), (\ref{second_Knuth}), (\ref{third_Knuth}), (\ref{fourth_Knuth}) and (\ref{fifth_Knuth})
as the Knuth-Taskin relations. 
Note that they 
are read either from the left to the right or conversely.
The main results Theorem 3.4 and Theorem 3.5 of [T]
amount in the $ r=0 $ case to the following:
 \begin{theorem}{\label{taskin}}
Suppose $ w, z \in W_n $. Then they belong to the same plactic class if and only if 
there is a sequence $ w_1, w_2 , \ldots , w_k \in W_n $ such that $ w = w_1, \, z = w_k $ and $ w_i 
\stackrel{KT}{\sim} w_{i+1} $ 
for $ i= 1 ,2 \ldots, k-1 $. In other words, the plactic classes are generated by 	
the Knuth-Taskin relations. 
\end{theorem}

The dual Knuth-Taskin relations are defined by $ w \stackrel{DKT}{\sim} w_1  $ if 
$ w^{-1} \stackrel{KT}{\sim} w_1^{-1}  $. If $ w $ and $ w_1 $ are written in word 
form, they do not act on neighboring elements, and as a matter of fact, they do not admit as simple 
a description as in the symmetric group case. On the other hand, since $ Q(w) = P(w^{-1}) $, 
we get an obvious dual version of 
the previous Theorem: 

\begin{theorem}{\label{dual_taskin}}
Suppose $ w, z \in W_n $. Then they belong to the same coplactic class if and only if 
there is a sequence $ w_1 , w_2, \ldots , w_k \in W_n $ such that $ w = w_1, \, z = w_k $ and $ w_i 
\stackrel{DKT}{\sim} w_{i+1} $ 
for $ i=1 ,2 \ldots, k-1 $.
\end{theorem}

\medskip
We need the following Lemma. 
\begin{lemma}{\label{W_b stable}} $W_b$ is stable under the Knuth-Taskin relations (\ref{first_Knuth}), (\ref{second_Knuth}), 
(\ref{third_Knuth}), (\ref{fourth_Knuth}) and (\ref{fifth_Knuth}).
\end{lemma}
\begin{proof} Assume that $ w \in W_b $ and write it in the form
$$ w = \underline{i}_1 \, \overline{a}_1 \, \underline{i}_2 \, \overline{a}_2 \, \ldots 
\underline{i}_l \, \overline{a}_l \, w_1 $$
where $w_1, \underline{i}_j $ are words, possibly 
empty, over $ I_n^+ $ for $ j=1,2, \, \ldots , \,l $ and
$ a_j > 0 $ for $ j=1,2, \, \ldots , \,l $. Note that we allow $ l = 0 $ corresponding to 
$ w = w_1 $. 
Write 
$$   \underline{i}_1 \underline{i}_2 \ldots \underline{i}_l w_1= i_1 i_2 \ldots i_k $$

Assume now that the Knuth-Taskin relation (\ref{first_Knuth})
acts in the 
$$  \underline{i}_1 \, \overline{a}_1 \, \underline{i}_2 \, \overline{a}_2 \, \ldots 
\underline{i}_l \, \overline{a}_l  $$
part of $w$.
We know from Theorem \ref{first_thm} that all $ \underline{i}_j  $ 
are increasing sequences over $ I_n^+ $ and that 
\begin{equation}{\label{we_know}} a_l < a_{l-1 } < a_{l-2} < \ldots < a_1 < 
\underline{i}_1 < \underline{i}_2 < \ldots < \underline{i}_l 
\end{equation}
where the inequalities hold for all elements of the subsequences, 
and so the pattern $ f(2) \, f(3) \, f(1) $ can only occur if $ f(1) = \overline{a}_r $ for some
$ 1 \leq r \leq l $ and $ f(3) = i_s $ for some $ s $. But then clearly (\ref{first_Knuth}) takes 
$ w$ 
to another element of $ W_b$.
Likewise we see that (\ref{first_Knuth}) acting in the pattern $ f(2) \, f(1) \, f(3) $
of $  \underline{i}_1 \, \overline{a}_1 \, \underline{i}_2 \, \overline{a}_2 \, \ldots 
\underline{i}_l \, \overline{a}_l  $
takes $w $ to another element of $ W_b $. 

\medskip
In the case of 
the Knuth-Taskin relation (\ref{second_Knuth}) acting in 
$$  \underline{i}_1 \, \overline{a}_1 \, \underline{i}_2 \, \overline{a}_2 \, \ldots 
\underline{i}_l \, \overline{a}_l  $$
we argue similarly. 
By the inequalities (\ref{we_know}), the only decreasing subsequences of  
$  \underline{i}_1 \, \overline{a}_1 \, \underline{i}_2 \, \overline{a}_2 \, \ldots 
\underline{i}_l \, \overline{a}_l  $ are of the form $i_r \, \overline{a}_s $
for some $ r , s $ and so in the pattern $ f(1) \, f(3) \, f(2) $ we have that $ f(3 ) = i_r $ for some $r$ whereas
$ f(2) =  \overline{a}_s $ for some $ s$. 
But since $ f(1) $ is less than $ f(2)  $ it must be $ \overline{a}_t $ for some $ t$ 
and so changing $ f(1) \, f(3) \, f(2) $ to $ f(3) \, f(1) \, f(2) $
gives another element of $ W_b$.
We argue similarly 
in case of the pattern $ f(3) \, f(1) \, f(2) $.

\medskip
Assume now that one of the Knuth-Taskin relation 
(\ref{first_Knuth}) or (\ref{second_Knuth})
acts in the $ w_1 $ part of $ w $. 
By the theory of the usual Robinson-Schensted 
algorithm, 
the length of the longest decreasing 
subsequence is preserved when the action is on words over $ I_n^+ $, 
and hence we get from Theorem 
{\ref{first_thm}} 
that (\ref{first_Knuth}) and (\ref{second_Knuth}) map $ w$ to an element of $ W_b $ in this case.

\medskip
We then 
consider 
the case where the action of one of the 
Knuth-Taskin relations
(\ref{first_Knuth}) and (\ref{second_Knuth})
involves both 
$ \underline{i}_1 \, \overline{a}_1 \, \underline{i}_2 \, \overline{a}_2 \, \ldots 
\underline{i}_l \, \overline{a}_l  $ and 
$ w_1 $. 
In that case 
$ \overline{a}_l $ must occur in first or second position of the relation. 

\medskip \noindent
{\bf Case $ f(2)  f(3)  f(1) $}: This case does not occur since $ f(1) $ would belong to $ w_1 $ and
would be less than $ \overline{a}_l $, which contradicts the fact that $ w_1 $ is a word over $ I_n^+$.

\medskip \noindent
{\bf Case $ f(2) f(1) f(3) $}: Using once more that the only decreasing 
subsequences of $ \underline{i}_1 \, \overline{a}_1 \, \underline{i}_2 \, \overline{a}_2 \, \ldots 
\underline{i}_l \, \overline{a}_l  $ are of the form 
$i_r \, \overline{a}_s $,
we get in this case that $ f(1) = \overline{a}_l $ 
whereas $ f(2 ) $ is unbarred. Applying the Knuth-Taskin relation (\ref{first_Knuth}) yields $ f(2) f(3) f(1) $, 
and hence $ \underline{i}_l $ changes to $ \underline{i}_l \,f(3) $, which is still increasing. 

\medskip \noindent
{\bf Case $ f(1) f(3) f(2) $}: In this case we have that $ f(1) = \overline{a}_l $ and $ f(3) $ and $ f(2) $ are 
unbarred, since $ f(2) \in w_1 $ and $ f(3) > f(2)  $. Thus also $ f(3) \in w_1 $.
The application of the Knuth-Taskin relation (\ref{second_Knuth}) changes $ f(1) f(3) f(2) $ to 
$ f(3) f(1) f(2) $ and hence $ \underline{i}_l $ changes to $ \underline{i}_l \, f(3) $. But no element of
$ \underline{i}_l $ can be bigger than $ f(3) $ for if $ i_r $ were such an element than we may 
assume it is the last one of $ \underline{i}_l $ and 
$ i_r \, f(3) \, f(2) $ would 
be a decreasing subsequence longer than three, inside 
$ a_l \, a_{l-1} \, a_{l-2} \, \ldots \, a_1 \, \underline{i}_1 \, \underline{i}_2 \,  \ldots 
\underline{i}_l \, w_1 $. Thus $ \underline{i}_l \, f(3) $ is increasing and we are done in this case 
as well.

\medskip \noindent
{\bf Case $ f(3) f(1) f(2) $}: We have $ f(1) = \overline{a}_l $. Using the Knuth-Taskin relation
(\ref{second_Knuth}) we have that $ f(3) f(1) f(2) $ changes to 
$ f(1) f(3) f(2) $ and thus $  \underline{i}_l  $ changes to $  \underline{i}_l \setminus f(3) $
which is clearly increasing. 

\medskip 
We next check that also the 
third Knuth-Taskin relation (\ref{third_Knuth}) takes $ w \in W_b $ to an element of $ W_b$. 
If $ w = w_1 $ then $ i_1 > i_2  $ since we are supposing that (\ref{third_Knuth}) 
acts in $ w $. Using Theorem \ref{first_thm} we then find that the image of $ w $ 
under (\ref{third_Knuth}), namely 
$ \overline{i_1} i_2  \ldots $, also belongs to $ W_b$. 

\medskip
In the remaining cases at least one of the two first elements of $ w $ must be negative.
Using Theorem \ref{first_thm} they are either of 
the form $ i_1 \, \overline{a_2} $ with $ i_1 > a_2 $ 
or $ \overline{a}_1 \, \overline{a}_2 $ with $ a_1 > a_2 $. 
But then from Theorem \ref{first_thm} once again we find in each case that the image of $ w$ under 
(\ref{third_Knuth}) also belongs to $ W_b $. 
\medskip

We finally show that the Knuth-Taskin relations 
(\ref{fourth_Knuth}) and (\ref{fifth_Knuth}) map $ w$ to an element of $ W_b $.
For this we assume that $ w \in W_b $ is either of the form $ w= \alpha $ or $ w= \alpha^{\prime} $ 
in the description of (\ref{fourth_Knuth}) and (\ref{fifth_Knuth}), and we let $ u, u^{\prime}$ and $ k $ 
be chosen correspondingly. 

\medskip
Let us first consider the Knuth-Taskin relation (\ref{fourth_Knuth}). 
We claim that if $ k \geq 2  $ and either $ \alpha  $ or $ \alpha^{\prime} $ belongs to $ W_b $, 
then the Knuth-Taskin relation (\ref{fourth_Knuth}) does not apply. 
Let us first verify this for $ k \geq 3 $. In that case 
$ u $ and $ u^{\prime} $ 
are both of the form 
\begin{equation}{\label{of_the_form}} a_{1,1} b_{1,1} a_{2,2} a_{2,1} b_{2,2} b_{2,1}  \ldots 
\end{equation}
where the dotted elements contain both $ a_{3,1} $ and $ b_{3,1} $ and hence both positive and negative numbers.
By the conditions on (\ref{fourth_Knuth}) 
we have the inequalities 
$ | b_{2,1} | <   | b_{2, 2} |  $, $ | a_{2,1} | <   | a_{2, 2} |  $ and $ | b_{1,1} | <   | a_{2, 2} |   <   | b_{2,2} |  $, 
among others.
If $ a_{i,j} > 0 $ (or equivalently $ b_{i,j} < 0 $) 
we get a contradiction with 
Theorem
{\ref{first_thm}} that implies $   | a_{2,2} |  < | a_{2,1}|  $ if $ \alpha $ or $ \alpha^{\prime} $ belongs to $W_b$.
If $ a_{i,j} < 0 $ (or $ b_{i,j} > 0 $) 
we also get a contradiction with 
Theorem
{\ref{first_thm}} that implies 	
$  | b_{2,2}  | <  | b_{2,1} | $, and the claim is proved for $ k \geq 3$.

\medskip 
In the case $ k = 2 $ we have 
$$u =  a_{1,1} b_{1,1} a_{2,2} a_{2,1} b_{2,2} b_{2,1} a_{3,2} a_{3,1} z
\,  \mbox{ and } \, 
u^{\prime} =  a_{1,1} b_{1,1} (-b_{2,2})  a_{2,2} a_{2,1} b_{2,1} a_{3,2} a_{3,1} z 
$$
and the same inequalities hold, that is 
$ | b_{2,1} | <   | b_{2, 2} |  $, $ | a_{2,1} | <   | a_{2, 2} |  $ and $ | b_{1,1} | <   | a_{2, 2} |   <   | b_{2,2} |  $.
We can then argue as above to show that $ \alpha $ cannot belong to $W_b$, indeed 
$ b_{i,j} < 0 $ implies by 
Theorem {\ref{first_thm}} that 
$ | a_{2,2} | <   | a_{2, 1} |  $ and $ b_{i,j} > 0 $ implies $ | b_{2,2} | <   | b_{2, 1} |  $.
Similarly, if $ \alpha^{\prime} \in W_b $ 
and $ b_{i,j} < 0 $ we get from Theorem {\ref{first_thm}} that $ | a_{2,2} | <   | a_{2, 1} |  $, which is a contradiction
and if $ b_{i,j} > 0 $ we get from Theorem {\ref{first_thm}} that $ | b_{2,2} | <   | a_{1,1} |  $,
which is also a contradiction.

\medskip
The only remaining possibility is $ k = 1 $. In that case we have 
$$ \begin{array}{lr}
u = a_{1,1} b_{1,1} a_{2,1} z,  & 
u^{\prime} = (- b_{1,1}) a_{1,1}  a_{2,1} z
\end{array}
$$
where $  | a_{1,1} | <   | a_{2,1} |  $.
If $ b_{1,1} < 0 $ we get by Theorem \ref{first_thm} that $ \alpha \in W_b $ iff $ \alpha^{\prime} \in W_b $.
On the other hand, if $ b_{1,1} > 0 $ we have that $ a_{1,1} $ and $  a_{2,1} $ are negative and hence 
by the inequality neither $ \alpha$ nor 
$ \alpha^{\prime} $ belongs to $W_b $, using Theorem \ref{first_thm} once again.
\medskip

We then finally treat the Knuth-Taskin relation (\ref{fifth_Knuth}). We proceed in the same way as for (\ref{fourth_Knuth}). 
If $ k \geq 3 $
we have that 
$ u $ and $ u^{\prime} $ 
both are of the form 
({\ref{of_the_form}})
where the dotted elements contain both positive and negative numbers
and where $ | b_{2,1} | <   | b_{2, 2} |  $, $ | a_{2,1} | <   | a_{2, 2} |  $ 
and $ | a_{1,1} | <   | b_{1,1 } |   <   | a_{2,2} |  $, among others.
The first two inequalities are the same as for (\ref{fourth_Knuth}) and so the argument used 
for (\ref{fourth_Knuth}) 
shows that neither $ \alpha $ nor $ \alpha^{\prime} $ can belong to $W_b $.

\medskip
When $ k = 2 $, we also use essentially the same argument as for (\ref{fourth_Knuth}). Indeed, we have 
$$ \begin{array}{l}
u = a_{1,1} b_{1,1} a_{2,2} a_{2,1} b_{2,2} b_{2,1} a_{3,3} a_{3,2} a_{3,1} b_{3,2} b_{3,1} z, \\
u^{\prime} =  a_{1,1} b_{1,1} a_{2,2} a_{2,1}( -a_{3,3}) b_{2,2} b_{2,1}  a_{3,2} a_{3,1} b_{3,2} b_{3,1} z
\end{array}
$$
and the same inequalities hold, that is 
$ | b_{2,1} | <   | b_{2, 2} |  $, $ | a_{2,1} | <   | a_{2, 2} |  $ and $ | a_{1,1} | <   | b_{1,1 } |   <   | a_{2,2} |  $. 
These inequalities ensure, using Theorem \ref{first_thm} as before, that $ \alpha $ does not belong 
to $ W_b$. On the other hand, if $ b_{i,j} < 0 $ and $ \alpha^{\prime} \in W_b $, we get by Theorem \ref{first_thm} that 
$ | a_{2,2} | <   | a_{2, 1} |  $, which is a contradiction. Finally, if $ b_{i,j} > 0 $ and 
$ \alpha^{\prime} \in W_b $ we get that $ | a_{2,2} | <   | b_{1, 1} |  $, which is also a contradiction.

\medskip
The only remaining case is now $ k = 1 $ where we have 
$$ \begin{array}{lr}
u = a_{1,1} b_{1,1} a_{2,2} a_{2,1} b_{2,1} z,  & 
u^{\prime} =    a_{1,1}(-a_{2,2}) b_{1,1}  a_{2,1} b_{2,1}  z
\end{array}
$$
where $ | a_{2,1} | <   | a_{2, 2} |  $ and $ | a_{1,1} | <   | b_{1,1 } |   <   | a_{2,2} |  $. 
These inequalities imply by Theorem \ref{first_thm} that neither $ \alpha  $ nor $ \alpha^{\prime} $ is in $W_b $.
The Lemma is proved.
\end{proof}

\begin{Cor}
$ W_b$ is a union of plactic classes and also a union of coplactic classes.
\end{Cor}
\begin{proof}
The previous Lemma amounts to saying that $ W_b $ is a union of plactic classes. 
But $ Q(w) = P(w^{-1} ) $ and $ W_b $ is stable with respect to $ w \mapsto w^{-1} $, hence 
$ W_b $ is also a union of coplactic classes.
\end{proof}

\medskip
For $ w \in W_n $ 
we define $ Sh(w) $ by $ Sh(P(w)) $ or, equivalently, 
by $ Sh(Q(w)) $.
Define 
$$ STD^{ \,\, 2}_{ \le 2 } := \{ (S,T) \in STD^2 \, | \, Sh(S) \mbox{ has less than two lines } \,\}. $$
We are now in position to prove the main Theorem of this section.  
\begin{theorem}{\label{first-main}}
Suppose that $ w \in W_n $. Then $ w \in W_b $ if and only if $ Sh(w) $ is a Young diagram 
of at most two lines. In other words, $ W_b $ is in correspondence with $ STD^{ \, \, 2}_{ \le 2 } $ 
under domino insertion.
\end{theorem}
\begin{proof}
Assume first that $ Sh(w) $ has at most two lines. Using Theorems
\ref{taskin} and \ref{dual_taskin} there is $ w_1 \in W_n $ related to $ w $ 
through a series of Knuth-Taskin or dual Knuth-Taskin relations such that $ P(w_1 ) $ and $ Q(w_1 ) $ both 
have one of the forms given in Figure \ref{two}
\begin{figure}
\includegraphics{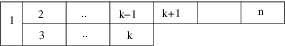}  \, \, \, \, \, \           \includegraphics{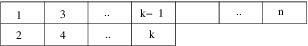}
\caption{}
\label{two}
\end{figure}
depending on the parity of the first line of $ Sh(w)$.
Under the domino insertion algorithm, the first tableau corresponds to 
$$ \underline{1} \, 3 \, 2 \, 5 \, 4 \, \ldots \, k \, k-1 \, k+1 \, k+2 \, k+3 \,  \ldots  \, n $$ 
whereas the second tableau corresponds to 
$$ 2 \, 1 \, 4 \, 3 \, 6 \, 5 \ldots \,k \, k-1 \, k+1 \, k+2 \, k+3 \ldots   \, n $$ 
Since they both belong to $ W_b $ we deduce from  
Lemma {\ref{W_b stable}} that $ w $ also belongs to $ W_b $ and one implication of the Theorem 
is proved.

\medskip
To prove the other implication we take $ w \in W_b $ and show that $ P(w) $ has at most two lines.
Write first $ w $ in the form 
$$ w = \underline{i}_1 \, \overline{a}_1 \, \underline{i}_2 \, \overline{a}_2 \, \ldots 
\underline{i}_u \, \overline{a}_u \, w_1 $$
where $w_1, \underline{i}_j $ are words over $ I_n^+ $ and
$ a_j > 0 $. We set  
$$ i_{1} \, i_{2} \, \, i_{3} \ldots \, i_{k} :=  \underline{i}_1 \, \underline{i}_2 \ldots \underline{i}_u.  $$
By Theorem \ref{first_thm} there is now a $ t $ such that 
$ P := P(\underline{i}_1 \, \overline{a}_1 \, \underline{i}_2 \, \overline{a}_2 \ldots 
\underline{i}_u \, \overline{a}_u ) $
is the domino tableau given in Figure \ref{three}.
\begin{figure}
\includegraphics{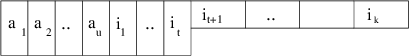}  
\caption{}
\label{three}
\end{figure}
Let $ w_1 = j_{1} \, j_2 \ldots \, j_{n-k-u} $ and let $ j_{i_1} \, j_{i_2} \ldots \, j_{i_r} $ 
be the subsequence of $ w_1 $ consisting of those elements $ j_i $ that are less then $ i_k $. Then by 
Theorem \ref{first_thm} we have that 
$ j_{i_1} \, j_{i_2} \ldots \, j_{i_r} $ is an increasing subsequence. Let $ j_{\iota_1} \, j_{\iota_2} \ldots
j_{\iota_s} $ be the subsequence of $ w_1 $ consisting of those elements that are positioned before 
$ j_{i_r} $ in $ w_1 $ and are bigger than $ i_k $. By 
Theorem \ref{first_thm} this is also an increasing subsequence. 
Setting $ z_1 := j_1 j_2 \ldots j_{i_r} $ and 
$ z_2 := j_{i_{r}+1} j_{i_{r}+2}    \ldots j_{n -k-u} $ we have obviously that $ w_1 = z_1 z_2 $. 
Moreover 
$ z_1 $ is a shuffle of its subsequences $ j_{i_1} \, j_{i_2} \ldots \, j_{i_r} $ and 
$ j_{\iota_1} \, j_{\iota_2} \ldots j_{\iota_s} $.
Let us first assume that this shuffle is trivial in the sense that 
$ z_1 = j_{i_1} \, j_{i_2} \ldots \, j_{i_r} j_{\iota_1} \, j_{\iota_2} \ldots j_{\iota_s} $.

\medskip
Let us consider the insertion of $ z_1 $ in 
$P$. If $ j_{i_1} $ must be 
entered in the two-line part of $ P $, say if $ a_1 < j_{i_1} < a_2 $, 
the resulting domino will be as in Figure \ref{four}, 
\begin{figure}
\includegraphics{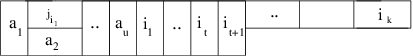}  
\caption{}
\label{four}
\end{figure}
that is, one vertical domino in $ P $ become horizontal, and the first horizontal domino 
becomes vertical. If $ j_{i_2} $ must also be entered in the two-line part of the tableau, the resulting 
tableau will look as in Figure \ref{five}
\begin{figure}
\includegraphics{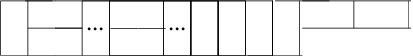}  
\caption{}
\label{five}
\end{figure}
where once again a vertical domino becomes horizontal and a horizontal becomes vertical. 
Since the sequence $ j_{i_1} \, j_{i_2} \ldots \, j_{i_r} $ is increasing this pattern is repeated 
until arriving at the elements that must be inserted in the one-line part of the tableau. These are 
inserted 
by bumping horizontal dominoes to the second line, giving tableaux of the form given in Figure \ref{six}.
\begin{figure}
\includegraphics{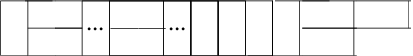}  
\caption{}
\label{six}
\end{figure}
We next describe the 
insertion of the other elements of $ z_1 $, those from $ j_{\iota_1} \, j_{\iota_2} \ldots
j_{\iota_s} $. But this is much simpler, since the element to be inserted will always be 
bigger than those so far inserted. It is therefore inserted   
as a horizontal domino at the end of the first line, without bumping. 

\medskip
This last description also shows 
that in general, when $ z_1 $ is a more complicated shuffle of 
$ j_{i_1} \, j_{i_2} \ldots \, j_{i_r} $ and 
$ j_{\iota_1} \, j_{\iota_2} \ldots
j_{\iota_s} $, 
the insertion of the elements of $ j_{\iota_1} \, j_{\iota_2} \ldots
j_{\iota_s} $, does not influence the insertion of the elements of
$ j_{i_1} \, j_{i_2} \ldots \, j_{i_r} $. 
We have thus proved that the insertion of all elements of $ z_1 $ gives a two-line domino 
tableau.

\medskip
Finally, we consider the insertion of the elements of $ z_2 $.
But the elements of $ z_2 $ are all bigger than 
the elements of $ P_1 = P(\underline{i}_1 \, \overline{a}_1 \, \underline{i}_2 \, \overline{a}_2 \ldots 
\underline{i}_u \, \overline{a}_u  z_1  )$
and so they are inserted as horizontal dominoes at the end 
of $ P_1 $. To be precise, the resulting domino tableau 
is simply the concatenation of the 
lines of $ P_1 $ and $ P(z_2) $. The Theorem is proved.
\end{proof}	

\medskip
In the remaining part of this section, we 
formulate a result which is a first strong indication 
of the connection between 
the empty core case of the [BGIL] conjectures and 
the representation theory of $b_n  $, 
where $b_n$ is the blob algebra mentioned in the introduction.

\medskip
Let 
$ \Par(n) $ denote the set of integer partitions of degree $ n $ with empty core and  
set $ \Par := \bigcup_{ n \geq 0 } \, \Par(n) $. Similarly, let 
$ \Bip(n) $ denote the set of bipartitions $ (\lambda, \mu ) $ of total degree $ n $ and 
set $ \Bip := \bigcup_{ n \geq 0 } \, \Bip(n) $. 
We denote by $ ST_{\emptyset}(n), ST_{\emptyset}, SBT(n) $ and $ SBT $ the set of standard (bi)tableaux with underlying 
shape in $ \Par(n), \Par, \Bip(n) $ and $ \Bip $.
For $ \lambda $ a partition we denote 
by $ {\mathcal Q}(\lambda ) $ 
the two-quotient of $ \lambda$, see for example [M] for a definition of it.
Then 
$ {\mathcal Q}(\lambda) \in \Bip(m) $ if $ \lambda \in \Par(2m) $ and $ \mathcal Q $ induces a bijection  
$$ {\mathcal Q} : \Par \rightarrow \Bip.$$
Following [BGIL] we define a partial order on 
$ \Bip $ by the rule 
$$ (\lambda, \mu ) \prec ( \tau, \nu  ) \, \, \, \, \mbox{iff} \, \, \, \,
{\mathcal Q}^{-1}(\lambda, \mu ) \lhd {\mathcal Q}^{-1}( \tau,  \nu ) $$
where $ \lhd $ refers to the usual dominance order on partitions.

\medskip
Let $ \Bip_1(n) $ denote the set of one-line bipartitions of total degree $ n $.
An element of $ \Bip_1(n) $ is of the form $ (\lambda, \mu ) =(( a), (n-a) )  $ for 
some positive integer $ a $ with $ 0 \le a \le n $.
We shall use the shorthand notation $ (a), (n-a ) $ for such $ (\lambda, \mu ) $
but reserve the notation $ (a, n-a ) $ for a conventional (two-line) partition.
Set $ \Bip_1 :=  \bigcup_{ n \geq 0 } \, \Bip_1(n) $.

\medskip
Define $ \Lambda_n := \{ -n, -n+2, \ldots, n-2 , n \} $. Then there is a bijection 
$$ f:\Bip_1(n) \rightarrow \Lambda_n, \, \, \, (a), (b) \mapsto a-b. $$
Let $ \prec $ (also) denote the order on $ \Lambda_n$ induced by $ f $, that is, for $ \lambda, \mu \in \Lambda_n $,
$ \lambda \prec \mu $ iff $ f^{-1}(\lambda ) \prec f^{-1}(\mu ) $. 

\medskip
Note that $ \Lambda_n $ is the parameterizing set for the quasi-hereditary category 
$ b_n \modc $ of $b_n$-modules.
The hereditary order is given by $ \lambda <_{qh} \mu $ iff $ | \lambda | > | \mu | $ for 
$ \lambda, \mu \in \Lambda_n $.	We now have the following result.
\begin{theorem}{\label{coideal}}
a) $ \Bip_1(n) $ is a coideal in $ \Bip $ with respect to $ \prec$. \newline
b) The order $ \prec $ on $ \Lambda_n $ is a refinement of $ <_{qh} $.	
\end{theorem}
\begin{proof}
In [CL] a bijection $ \overline{\mathcal Q}:  SDT \rightarrow SBT $ is described. It induces
$ {\mathcal Q}: \Par(2n) \rightarrow \Bip(n) $ by taking shapes. 
One then checks the following formulas
$$ \begin{array}{ll}
{\mathcal Q}^{-1}: (a),(b) \mapsto (2a, 2b )  & \, \, \mbox{ for } \, a \geq b \\
{\mathcal Q}^{-1}: (a),(b) \mapsto ( 2b-1, 2a +1) & \, \, \mbox{ for } \, a <  b
\end{array}
$$
We deduce that $ {\mathcal Q}^{-1} ( \Bip_1(n) ) $ consists of all partitions of $ 2 n $ of at most two lines
and thus $ \Bip_1(n) $ indeed is a coideal in $ \Bip $ with respect to $ \prec $ as claimed in $a)$.

\medskip \noindent
In order to prove $b)$ we note that the above formulas give 
$$  (n), (\emptyset)  \succ   (\emptyset) , (n)  \succ  (n-1), (1)  \succ (1), (n-1) 
\succ  (n-2), (2)\succ \ldots  
$$
The last term is $ (\frac{n}{2}), (\frac{n}{2}) $
or 
$ (\frac{n-1}{2}), (\frac{n+1}{2}) $ depending on the parity of $ n$.
The statement of $ b) $ follows from this. In fact we see that 
the only difference between 
$ \prec $ and $ <_{qh} $ is that $ -\lambda  \prec \lambda $ if $ \lambda \in \Lambda_n $ and 
$ \lambda > 0 $, whereas 
they are noncomparable with respect to $ <_{qh} $.

\end{proof}

\section{Cell theory in ${ {\mathcal H}_n }$}
The fundamental text on cell theory for Hecke algebras 
with unequal parameters is Lusztig's book [Lu2]. Since we are here interested in the special $ B_{n} $ 
case, we shall follow the notation used in [BGIL].
Let therefore $ \Gamma $ be a finitely generated free Abelian group containing the elements $ a, b$ and 
let $ < $ be a total order on $ \Gamma$, making it into an ordered group.
We use exponential notation for the elements of $\Gamma $, writing $ e^g $ for  
$ g \in \Gamma$. Define  
$ \q := e^a $ and $ \Q := e^b $. Let $ \A $ be the $ \Z $-algebra
$ \A :={\mathbb Z}[\Gamma] $. 
The Hecke algebra $ {\mathcal H}_n = {\mathcal H}_n(\Q,\q)  $ of type $ B $ is the 
$\A$-algebra generated by $ T_0, T_1, \ldots , T_{n-1} $ subject to the 
relations 
$$ \begin{array}{ll} 
T_i T_{i - 1 } T_i = T_{i - 1 }  T_{i  } T_{i - 1 } & \mbox{for } i = 2,3, \ldots , n-1  \\
T_0 T_1 T_0 T_1 = T_1 T_0 T_1 T_0 & \\
T_i T_j = T_j T_i \mbox{ for } | i -j | > 1 &\\
(T_i -\q ) (T_i +\q^{-1} ) = 0, & (T_0 -\Q ) (T_0 +\Q^{-1} ) = 0.
\end{array} 
$$
The frequently used ground ring in the literature $\Z[ \Q,\Q^{-1},  \q, \q^{-1} ] $ 
is obtained as a special case of the above by setting $ \Gamma := \Z a \oplus \Z b $. 
The Hecke algebra 
defined over this ground ring is 
called the generic Hecke algebra.

\medskip

Assume that $ f: \Gamma \rightarrow \C^{\times} $ is a group homomorphism. Then 
$ f $ extends canonically to an algebra homomorphism $ f: \A \rightarrow  \C $ and we can define 
the specialized Hecke algebra
$ {\mathcal H}_{n,f} := {\mathcal H}_{n} \otimes_{ \A }{ \C} $.
For example $ f(g) = 1, \, \forall g $ gives the group algebra $ {\mathcal H}_{n,f} = \C W_n $.

\medskip
Define elements $ C_i $ of $ {\mathcal H}_n $ by $ C_0 := T_0 - \Q $ and $ C_i := T_i - \q $
for $ i = 1, 2, \ldots, n-1 $.
Let $ J_n $ be the following ideal of $ {\mathcal H}_n $
$$ J_n := \langle \, C_1 C_2 C_1 - C_1, \, \,  
 C_1 C_0 C_1 - [2]_{ \frac{\Q}{\q}}C_1 \, \rangle $$ 
where $ [n]_x := x^{n-1} + x^{n-3} + \ldots + x^{-n+3} +
x^{-n+1} $ is the usual Gaussian integer.  
We then define the Temperley-Lieb algebra of type $ B $
as $$ TLB_n := {\mathcal H}_n / J_n.$$
In the case of the generic Hecke algebra, this definition already appears in [GL1]
where $ TLB_n $ is also referred to as the blob algebra, but actually 
it differs slightly from the presentation of the blob algebra $ b_n $ that is used in eg. 
[MR] and [RH]. 
Let us be more specific about the relationship. 

\medskip
Let $ k $ be a field and
assume that $ q \in k^{\times}, q \not=1, -1  $ and $ m \in \Z$.
In [RH] and other references $b_n = b_n(q,m) $ is defined as the 
$ k $-algebra on generators $ U_0, U_1, U_2 \, \ldots, U_{n-1} $
and relations 
$$ \begin{array}{l}
U_i U_{i + 1 } U_i = U_i \,\, \,   \mbox{ for } i = 1 , 2 \ldots , n-2 \\
 U_{i + 1 } U_i U_{i + 1 } = U_{i+1} \,\, \, \mbox{ for } i = 1 , 2 \ldots , n-2 \\
U_1 U_{0 } U_1 = [m-1 ] U_1 \\ U_i^2 = -[2] U_i \mbox{ for } i = 1 , 2 \ldots , n-1 \\
U_0^2 = -[m] U_0, \, \, \,   U_i U_j  = U_j U_i  \mbox{ for } | i-j | > 1 
\end{array}
$$ 
where $ [a] = \frac{ q^{a} - q^{-a}}{q - q^{-1} } $. 
The following Lemma relates this to $ TLB_n $.
\begin{lemma}{\label{relates}}
Suppose 
$ k := \C $. Assume $ q \in \C^{\times } \setminus \{ 1,-1 \} $ and  
set $ Q := i q^m $. Define $ TLB_{n, q, Q} := TLB_n \otimes_{\A} \C $ where 
$ \C $ is made into an $ \A $-algebra via $ f : \Gamma \rightarrow \C^{\times} $ such 
that $ f(a) = q, f(b) = Q $. Then the rules 
$$ 
C_i \mapsto U_i,  \, \, \, \,   i = 1 , 2,  \ldots,  n-1, \, \, \, \,
C_0 \mapsto  i (q-q^{-1}) \,   U_0  
$$ 
define an isomorphism $ g: TLB_{n, q, Q} \rightarrow b_n(q, m)  $.
\end{lemma}
\begin{proof}
It is shown in 
Proposition 2.1 of [CGM] that $ TLB_n = { \mathcal H}_n / J_n $ is free
over $ {\mathcal A} ={\mathbb Z} [ \Q, \Q^{-1},  \q, \q^{-1} ] $ and hence 
the proof is only a matter of checking the relations. 
\end{proof}

\medskip
For $ w \in W_n $ we define $ T_w := T_{i_1 } T_{i_2 } \ldots T_{i_N } $ where 
$ w = s_{i_1} s_{i_2} \ldots s_{i_N} $ is a reduced expression. By the relations, $ T_w $ 
is independent of the reduced expression. Moreover, $ T_w $ is invertible since 
$ T_i $ is invertible for all $ i $; indeed we have 
\begin{equation}{\label{inv}}
 \begin{array}{ll} 
T_0^{-1} = T_0 -\Q + \Q^{-1}, &  
T_i^{-1} = T_i -\q + \q^{-1} \,  \mbox{ for } i = 1, 2, \ldots, n-1.
\end{array} 
\end{equation}
The bar involution $ h \mapsto \overline{h} $ on $ {\mathcal H}_n $ is the ring automorphism given by
$$ \begin{array}{lll} 
T_w \mapsto T_{w^{-1}}^{-1}, & \q \mapsto \q^{-1},      &  \Q \mapsto \Q^{-1}.  
\end{array} 
$$

\medskip
Recall now that $ \Gamma $ is endowed with a total order $ < $. 
Using it Lusztig introduces in [Lu1, Lu2] a Kazhdan-Lusztig type basis $ C_w, \, w \in W_n $ for $ {\mathcal H}_n $.
It is uniquely defined by the conditions
$$ 
\begin{array}{ll} 
\overline{C_w} = C_w, &  
C_w - T_w \in \bigoplus_{w \in W_n}  \A_{ >  0} T_w
\end{array} 
$$
where $ \A_{>  0} := \sum_{\gamma \in \Gamma, \gamma > 0} \C \,e^{\gamma} $. 

\medskip
Associated with the basis $ C_w $ there is a preorder $ \le_L  $ on $ W_n $, 
generated by $ y \le_L w $ if $ C_y $ appears in the expansion of $ C_{s_i } C_y $ in the $ C_w $-basis. 
The associated equivalence relation is denoted $ \sim_L $ and its classes left cells. 
Thus, $ z \sim_L w $ if $ z \le w $ and $ w \le z $.
Similarly we 
define the preorders $  \le_R  $ and $ \le_{LR} $ and the equivalence relations $ \sim_R $ and $ \sim_{LR} $.
The associated classes are called right cells and two-sided cells.

\medskip
We shall always assume that $ a $ and $ b $ are positive in $ \Gamma $ and so we get 
by the equations (\ref{inv}) the following formulas
$$ \begin{array}{ll}
C_{s_0} = T_{0} -\Q, & C_{s_i} = T_i -\q \, \, \, \mbox{for } i = 1,2, \ldots, n-1.
\end{array}
$$
In other words, we have that $  C_{s_i} = C_{i}  $. 

\medskip
Assume that $ b \not\in \{ a, 2a, \ldots , (n-1) a \} $. Let $ r \in { \mathbb N }
\cup \{  0 \} \cup \{  \infty \}$
be given by $ r a < b < (r+1) a $ or 
$ r := \infty $ if $ b>  (n-1)a $. 
According to the conjectures in [BGIL], the representation theory of $ {\mathcal H}_n $ 
should only depend on $ \Gamma, a $ and $ b$ through $r$.

\medskip
Let us consider the following $\mathcal A$-submodule of $ {\mathcal H}_{n} $
$$ {\mathcal J}_n := \Span_{\mathcal A}  \{ \, C_w \, | \,  w \not\in W_b \,  \}. $$
The next Theorem is the main result of this section. In order to formulate it, we recall 
that $ c^+ $ of Conjecture A of [BGIL] is the statement that 
$$ y \leq_{ \mathcal L R} w \Longleftrightarrow Sh(y) \leq Sh(w). $$

\begin{theorem}
Assume that $ r = 0 $ and 
assume that part $ {\rm c}^+ $ of Conjecture A of [BGIL] is valid for $ r = 0$. Then 
we have that 
$ {\mathcal J}_n = J_{n} $.
\end{theorem}
\begin{proof}
Since $ {\rm c}^+ $ is assumed to be true 
we have that $ \le_{\mathcal LR} $ is given by dominance order under domino insertion. 
Combining with Theorem {\ref{first-main}} we get that 
$ {\mathcal J}_n $ is an ideal in $ {\mathcal H}_n $.

\medskip
In order to show that $ J_n \subset {\mathcal J}_n $ it is then enough to verify that the generators of $ J_n $ belong 
to $ {\mathcal J}_n $. Now we have 
$$ \begin{array}{r}
C_1 C_2 C_1  = (T_1 -q )  (T_2 -q )  (T_1 -q )  = \\ T_1 T_2 T_1 - q \,T_1 T_2 - q \,T_2 T_1 
+ q^{2} \,T_1 + q^{2}\, T_2 - q^{3} +  T_1 -q 
\end{array}
$$
and hence 
$$ C_{s_1 s_2 s_1 } = C_1 C_2 C_1 - C_1 $$
On the other hand, $ P(s_1 s_2 s_1) $ has the form given in Figure \ref{seven}
\begin{figure}
\includegraphics{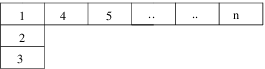}  
\caption{}
\label{seven}
\end{figure}
and so $ C_{s_1 s_2 s_1 } \in {\mathcal J}_n $. 
Similarly we have 
$$ 
\begin{array}{r}
C_1 C_0 C_1 = (T_1 -q) (T_0 -Q) (T_1 -q) = \\ 
T_1 T_0 T_1 - q \,T_0 T_1  +  Q q^{-1} \, T_1 -Q
- q \,T_1 T_0  + q^2 \,T_0  + q Q \,  T_1    -q^2 Q
\end{array}
$$
But $ -a +b < 0 $ and so $ Qq^{-1} \not\in{\mathcal A}_{> 0} $ and we must subtract $ [2]_{\frac{Q}{q}} \, C_1 $ to get 
$ C_{s_1 s_0 s_1 } $. Hence 
$$ C_{s_1 s_0 s_1 } = C_1 C_0 C_1 - [2]_{\frac{Q}{q}} \, C_1. $$
But $ P(s_1 s_0 s_1) $ is as in Figure \ref{eight}
\begin{figure}
\includegraphics{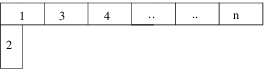}  
\caption{}
\label{eight}
\end{figure}
and so indeed $ C_{s_1 s_0 s_1 } \in {\mathcal J}_n $. 

\medskip
Let $ K $ be the kernel of the projection map  $ \pi :{\mathcal H}_{n} / J_n 
\rightarrow {\mathcal H}_{n} / {\mathcal J}_n $. 
We need to show that $ K = 0 $. 
Since $ \pi $ is surjective, it is enough to prove that 
${\mathcal H}_{n} / J_n $ and ${\mathcal H}_{n} / {\mathcal J}_n $ are free over $ \mathcal A $ of the same rank. 

\medskip 
As mentioned above, 
$ TLB_n = { \mathcal H}_n / J_n $ was shown in [CGM] to be free over 
the ground ring ${\mathbb Z} [ \Q, \Q^{-1},  \q, \q^{-1} ] $.
Its rank is given by the cardinality of 
the diagram basis and can also be read off from the 
Bratelli diagram for $ TLB_n $. It is $$ \rank  {\mathcal H}_{n} / J_n    = 
\sum_{ i=0}^{ n} \binom{ n}{i}^2. $$ 
On the other hand, since $ \{  C_w \} $ is a basis of $ {\mathcal H}_n $ we have that  
$  {\mathcal H}_{n} / {\mathcal J}_n $ is free over $ \mathcal A $ with rank 
$$ \rank {\mathcal H}_{n} / {\mathcal J}_n  = 
 |  W_b  |. $$
Recall the bijection $ \overline{\mathcal Q}: SDT \rightarrow SBT$ from [CL].
By the proof of Theorem {\ref{coideal}}, 
it restricts to a bijection between  
standard domino tableaux in $ STD(n) $ with less than two lines and 
one-line standard bitableaux with shape in $ \Bip_1(n) $.
The number of pairs of one-line 
bitableaux of shape $ (i, n-i ) $ is 
$ \binom{n}{i}^2 $
and so we conclude that $ \rank {\mathcal H}_{n} / { J}_n = \rank {\mathcal H}_{n} / {\mathcal J}_n $, as needed.
\end{proof}
\begin{remark} 
It is useful to observe that for the above proof to work, actually 
only '$ \Longrightarrow $' of part $ c^+ $ of Conjecture A in [BGIL] is needed.
\end{remark} 
\begin{Cor}
Assume that $ \Gamma = \mathbb Z $ with the standard order and that $b = 1 $ and $ a= 2 $. Then 
$ {\mathcal J}_n = J_{n} $.
\end{Cor}
\begin{proof}
By Remark 4.1 of [BJ], which on the other hand relies on [Lu1], we get that $ c^+ $ of Conjecture A 
of [BGIL] is valid under the assumptions.
We then apply the Theorem.
\end{proof}
In order to apply the Corollary, 
we shall from now on assume that $ \Gamma := \Z $ with the standard order, 
and that $ b:= 1 $, $ a = 2 $. Although this does not cover all of 
the $ r = 0 $ case of [BGIL] we shall, somewhat misleadingly, refer to it that way. 

\medskip
We need both versions of the blob algebra. 
Hence, in order for Lemma \ref{relates} and the Corollary to work we  
impose the following conditions on $q, Q  $ 
\begin{equation}{\label{condition}}
q \mbox{ \it{is a primitive {\it l}'th root of unity}}, \, l > 2, \,
Q := i q^m, \,  q= - q^{2m}.
\end{equation}
Note that the conditions imply that $ l $ is even.
They will be satisfied for example if
$ l = 2 (2m -1 ) $.

\medskip
{\it We choose from now on $ q, Q, m , l $ satisfying (\ref{condition}). We use the notation 
$ {\mathcal H}_{n, q, Q} $
for the specialized Hecke algebra $ {\mathcal H}_{f} $ with respect to these choices. Similarly, we write 
${\mathcal J}_{n, q, Q} $ for $ {\mathcal J}_f $ and ${ J}_{n, q, Q} $ for $ { J}_f $}. 

\begin{Cor}
We have $ TLB_{n, q, Q} =  {\mathcal H}_{n, q, Q}/ {\mathcal J}_{n, q, Q} = b_n (q, m) $.
\end{Cor}
\begin{proof}
This follows from the Theorem and Lemma \ref{relates}.
\end{proof}

\section{Representation theory}
In this section we use the results of the previous sections to study the representation theory 
of $ b_n$. Our main result is that the cell modules in the $ r=0 $ case are the standard modules 
for $ b_n$. 

\medskip
Recall that $ [2] \not= 0 $ so that we can 
define $ e= -\frac{1}{[2] }U_{n-1} $. This is an idempotent of $ b_n $ and we have 
that $ e b_n  e \cong b_{n-2} $. Hence it gives rise to the 
localization functor $$ F: b_n \modc \rightarrow b_{n-2} \modc, \, \, \, 
M \mapsto e M .$$
$F$ is exact, it has as left adjoint functor the globalization functor $ G $
$$ G: b_{n-2} \modc \rightarrow b_{n}  \modc, \, \, \, 
M \mapsto b_n e \otimes_{e b_n e}  M. $$
Recall that $ \Lambda_n := \{ -n, -n+2, \ldots, n-2 , n \} $ is the parameterizing set for 
the quasi-hereditary category $ b_n \modc $. Let $ \Delta_n( \lambda ) \in b_n \modc $ denote the standard 
module associated with $ \lambda \in \Lambda $. We have that 
\begin{equation}{\label{functors}}
\begin{array}{l}
F \Delta_n(\lambda) \cong \left\{ \begin{array}{ll} 
\Delta_{n-2}(\lambda) & \mbox{if} \,\,\,   \lambda \in \Lambda_{n}  \setminus \{\pm n\}\\
0  & \mbox{otherwise}  \end{array} \right. 
\\
G   \Delta_n(\lambda) \cong \, \, \,\, \,\,\,
\Delta_{n+2}(\lambda)  
\end{array}
\end{equation}
and $ \Delta_{n}(\pm n ) \cong L_n(\pm n)  $
where $ L_n(\lambda) $ is the simple module given by $ \lambda $.
This implies the universal property for $ \Delta_n(\lambda) $ as 
the projective cover of $ L_n(\lambda) $ in the truncated subcategory of $b_{n}  \modc$
consisting of modules with composition factors of the form $ L_n(\mu ) $ with $ \mu \leq \lambda $.

\medskip
Let now $ w_n \in W_{b} $ and denote by $ {\frak C}= {\frak C}_{w_n} \subseteq W_{b} $ its left cell. 
Consider the following ideals of $ {\mathcal H}_n $
$$ {\mathcal I}_{\leq_{\frak C} \, w_n} := \Span_{\C}   \{ C_w | w  \leq_L w_n  \}, \,  \,  \,  \,  
 {\mathcal I}_{<_{\frak C} \, w_n} := \Span_{\C}  \{ C_w | w  \le_L w_n, w \not\in {\frak C}  \} $$
and define the cell module $$ {\mathcal V}_{w_n} := {\mathcal I}_{\leq_{\frak  C} \, w_n} /  {\mathcal I}_{<_{\frak C} \, w_n} .$$
Since conjecture A of [BGIL] is true in the $ r=0 $ case, 
we get by the results of the previous section that 
$ {\mathcal V}_{w_n}  $ is a $ b_n$-module. A basis 
for $ {\mathcal V}_{w_n}  $ is given by the classes of $ C_w $ for $ w \in {\frak C} $.

\medskip
Recall from the previous sections that $W_n $ 
is realized as the subgroup of the symmetric group on the 
elements $ I_n $ generated by $ s_0 := (-1,1) $ and 
$ s_i := (i,i+1)(-i,-i-1) $. Let us denote by $ \iota $ the associated 
injection of groups $ \iota: W_n \rightarrow S_{I_n} =S_{2n}$:
$$ \iota(s_0) = (1, -1), \, \, \, \iota(s_i) = (i,i+1) (-i, -i+1). $$

According to the last Theorem of [Lu1] (on page 111), each left cell $ {\frak C} $ of $W_n $ is now of the form 
$ {\frak C} = \iota^{-1}(\widetilde{{\frak C }}) = \widetilde{{\frak C }}  \cap W_n  $ 
where  $ \widetilde{\frak C } $ is a left cell of $ S_{I_n} $; this relies heavily on $ r=0 $.

\medskip
By [KL,A1,G], the left cells on $S_{I_n} = S_{2n} $ can be described using the usual Robinson-Schensted correspondence
when we use the natural order on $ I_n $, given by
$$ \overline{n} < \ldots < \overline{2} < \overline{1} < 1 < 2 < \ldots  < n. $$
We need the following Lemma. 
\begin{lemma} Let $ {\frak C} $ be a left cell in $ W_n $. Assume that 
$ {\frak C} \subset W_b $ and that 
$ {\frak C} = \widetilde{\frak C }  \cap W_n  $ where $ \widetilde{\frak C } $ is a left cell of 
$ S_{I_n} $. Then under the Robinson-Schensted bijection on $ S_{I_n} $ with respect to 
the above order on $ I_n $, $ \widetilde{\frak C } $
corresponds to a tableau in at most two lines.
\end{lemma}
\begin{proof} 
Let $P^{\prime} $ and $Q^{\prime} $ denote the $ P $ and $ Q$-parts of the Robinson-Schensted 
correspondence on $ S_{2n} $. For $ z, z_1 \in \widetilde{\frak C } $ we have 
$Q^{\prime}(z ) = Q^{\prime}(z_1 ) $ and $P^{\prime}(z ) $ and $Q^{\prime}(z ) $ have the same shape.
Assume now that $ w \in W_b $ and write it in word form as $ w= i_{1} i_{2} \ldots i_{n} $
with $ i_{j} \in I_n $. 
We then have 
$$ \iota(w) = \overline{w}^{op}  w $$ 
where $ \overline{w}^{op}:= {\overline{i}}_n \,  {\overline{i}}_{n-1} \,  \ldots {\overline{i}}_{1} $
and so 
$ P^{\prime}(\iota(w)) = P^{\prime}(\overline{w}^{op}  w) $.

\medskip
We now appeal to the description of $ W_b $ given in Theorem 
{\ref{first_thm}}. Using it, there are no decreasing subsequences of $ \overline{w}^{op}  w $ 
of length three or more, and 
thus $ P^{\prime}(\overline{w}^{op}\,  w) $ has at most two lines. Indeed, consider the graph of Figure \ref{nine}.
\begin{figure}
\includegraphics{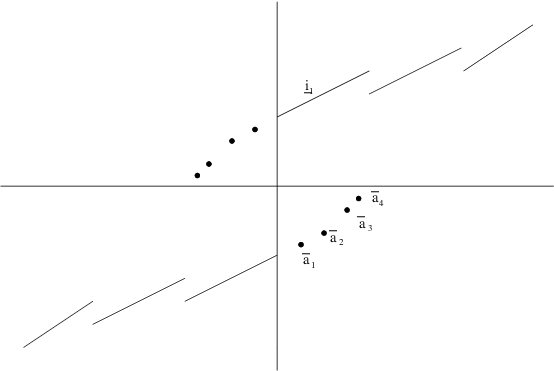}  
\caption{}
\label{nine}
\end{figure}
It represents $ \overline{w}^{op}  w $ in the case where 
$ \overline{a_1}, \ldots ,\overline{a_4} $ are the only negative numbers in $ w$, that is $ l= 4 $ in 
the notation of Theorem {\ref{first_thm}}. 
The restriction of the graph 
to the quadrants I and IV represents $ w $ and the restriction of the graph 
to the quadrants I and II represents $ w^l = w^4  $ in the notation of Theorem {\ref{first_thm}}.
For simplicity, the straight lines of quadrants I and III represent 
sequences of increasing numbers, where $ \underline{i}_1 := i_1 \ldots i_{k_4 -1} $.
From Theorem {\ref{first_thm}} we know that $ w^l $ has no decreasing subsequences of length 3 or more, 
and that these are all positioned after $ \underline{i}_1 $ in the graph. Hence 
the discontinuity points of quadrant I will increase as indicated. We now conclude that $ \overline{w}^{op}  w $
has no decreasing subsequences of length 3 or more, as claimed. The general case
is treated the same way.
\end{proof}

\begin{lemma}   
a) Assume that $ U_{n-1} {C}_{w_n} \not= 0 $. Then there exists $ w_{n-2} \in {\frak C}_{w_n} 
\cap W_{n-2} $ and a scalar $ a \in \C \setminus \{ 0 \} $ such that $ U_{n-1} C_{w_n} =  a \,U_{n-1} C_{w_{n-2}}$. \newline 
b) Assume $ U_{i} C_{w_n} \in {\mathcal V}_{w_n} \setminus \{0\}   $ for some $ i > 0 $. Then there
exists $ z \in {\frak C}_{w_n} $ and a scalar $ a \in \C \setminus \{ 0 \} $ such that $ U_{i} C_{w_n} = a C_z $. \newline 
c) Assume that $ U_{n-1} {\mathcal V}_{w_n} = 0 $. Then $ U_{i}\, {\mathcal V}_{w_n} = 0 $ for all $ i > 0$.
Moreover $ {\mathcal V}_{w_n} \simeq \Delta_n(\pm n ) $, specially
$ \dim {\mathcal V}_{w_n} =1$.
\end{lemma} 
\begin{proof}
Take $ w_n \in {\frak C}_{w_n} = {\frak C}  $ and let $ C_{w_n} \in {\mathcal H}_{n} $ be the associated 
Kazhdan-Lusztig element. Then we have 
\begin{equation}{\label{structure_constant1}} 
U_{n-1} C_{w_n} = C_{s_{n-1}} C_{w_{n}} = \sum_{z \in W_n}  N_{n-1, w_n, z} \,C_z 
\end{equation} 
where $ N_{n-1, w_n, z} $ are the structure constants in $ {\mathcal H}_n $ with respect to the $ C $-basis.
Let $ {\mathcal H}_{2n}  $
be the Hecke algebra associated to $ S_{2n} $, with parameter $ q $, and let us denote by $ \widetilde{C}_w $ the 
usual one-parameter Kazhdan-Lusztig element for $ w \in S_{2n}  $.
If $ w \in W_n $ we write
$ \widetilde{C}_w := \widetilde{C}_{\iota(w)}$. Then we have	
\begin{equation}{\label{structure_constant2}} 
\widetilde{C}_{s_{n-1}} \widetilde{C}_{w_n}= 
\sum_{z \in S_{2n}}  \widetilde{N}_{n-1, w_n, z} \,\widetilde{C}_z 
\end{equation} 
where $ \widetilde{N}_{n-1, y, z} $ are the structure constants in
$ {\mathcal H}_{2n} $ with respect to its $ \widetilde{C} $-basis.
Lusztig shows in this setting in [Lu1] that 
\begin{equation}{\label{lu}}
 \mbox{ if } z \in W_n \mbox{ and  } {N}_{n-1, w_n, z} \not= 0 \mbox{ then } 
\widetilde{N}_{n-1, w_n, z} \not=0. 
\end{equation}
\medskip
Now we have  
$$
\widetilde{C}_{s_{n-1}} = 
(T_{(n-1,n)} -q ) (T_{(-n+1,-n)} -q ) = U_{(n-1,n)} U_{(-n+1,-n)}.
$$
Reducing (\ref{structure_constant1}) modulo   
$ {\mathcal I}_{<_{\mathcal C} \, w_n} $ we get the corresponding equation in ${\mathcal V}_{w_n} $:
\begin{equation}{\label{structure_constant3}}	
U_{n-1} C_{w_n} = \sum_{z \in {\frak C}}  N_{n-1, w_n, z} \,C_z  \mbox{  modulo  } 
 {\mathcal I}_{<_{\mathcal C} \, w_n}. 
\end{equation}
But $ {\frak C} = \widetilde{\frak C} \cap W_n $ and so by (\ref{lu}) any $ z $ occurring in this sum with  
$N_{n-1, w_n, z} \not= 0$ gives a nonzero $\widetilde{N}_{n-1, w_n, z} $ in
\begin{equation}{\label{structure_constant4}} 
\widetilde{C}_{s_{n-1}} \widetilde{C}_{w_n}= U_{(n-1,n)} U_{(-n+1,-n)} \widetilde{C}_{w_n} =
\sum_{z \in S_{2n}}  \widetilde{N}_{n-1, w_n, z} \,\widetilde{C}_z 
\mbox{  modulo  } 
 {\frak I} 
\end{equation} 
where 
$$ {\frak I } := \Span_{\C} \{ \widetilde{C}_w | \, w \in S_{2n}, \, w \le_L \widetilde{\frak C}, \, 
w \not\in \widetilde{\frak C} \}. $$
But using the previous Lemma we may consider ({\ref{structure_constant4}}) as an equation in 
a cell module $ \Delta_{2n}(k) $ for the 
Temperley-Lieb algebra $TL_{2n} $.

\medskip
Let us now show a). We have $ N_{n-1, w_n, z} \not= 0 $ and so $ \widetilde{N}_{n-1, w_n, z} \not= 0 $.
But by [FG] we know that $ \widetilde{C}_z = U_{\iota(z)} \mbox{ modulo } {\frak I } $,
where as usual $ U_w := U_{i_1} \ldots U_{i_r} $ for $ w= s_{i_1} \ldots s_{i_r} $.
Using the diagram presentation of $ \Delta_{2n}(k) $ we now deduce that 
$ \iota(z) = s_{(-n+1,-n)} \, z_1 \, s_{(n-1,n)} $ where $ z_1 \in S_{I_{n-2}} $ and hence 
$$ z_1 = s_{(-n+1,-n)} \, \iota(z) \, s_{(n-1,n)}
=  \iota(z s_{n-1}) 
 \in \iota(W_{n}) \cap S_{I_{n-2}} =  \iota(W_{n-2})   $$ 
and a) is proved.

\medskip
We then show b). For each $ z $ with 
$  N_{i, w_n, z} \not= 0 $ we have by (\ref{lu}) that $  \widetilde{N}_{n-1, w_n, z} \not= 0$.
But using [FG] once more, at most one $ z $ can give 
$  \widetilde{N}_{n-1, w_n, z} \not= 0$, proving b).

\medskip
Let us then show c). By the previous sections, 
$ {\mathcal V}_{w_n} $ is a module for $ b_n $. 
Since $ F {\mathcal V}_{w_n}  = U_{n-1} {\mathcal V}_{w_n}  = 0 $, it follows from the general representation 
theory of $b_n $ that 
$$ {\mathcal V}_{w_n}  = \Delta_{n}(n)^k \oplus \Delta_{n}(-n)^l $$
for certain multiplicities $ k, l$. Since $ {\mathcal V}_{w_n} $ is a cell module, the products
$ C_{s_{i_1}}  \,\ldots \, C_{s_{i_k}} \, C_{w_n } $ generate $ {\mathcal V}_{w_n} $.
But by assumption only $ C_{s_0}^k \, C_{w_n } = U_0^k \,   C_{w_n } \in {\mathcal V}_{w_n} $
can be nonzero and since $ U_0^k $ is a scalar multiple of $ U_0 $ we conclude 
that $ k= 1, l=0 $ or $ k= 0, l=1 $ and so $ \dim  {\mathcal V}_{w_n} =1 $. The Lemma is proved. 
\end{proof}

We are now in position to prove our main Theorem.
\begin{theorem}{\label{main}}
Assume that $q$ is a primitive $l$'th root of unity such that $ q= - q^{2m} $ and 
$ Q := i q^m  $.
Let $ {\frak C} = {\frak C}_{w_n} $ be a left cell for $ W_n $ and let 
$ {\mathcal V }_{w_n} $ be the corresponding cell module. Then we have an isomorphism of $ b_n$-modules
$$  {\mathcal V }_{w_n} \simeq \Delta_n ( \lambda ) $$
where $ \lambda = a-b $ for $ {\mathcal Q}(Sh(w_n ))= (a),(b)$.
\end{theorem}
\begin{proof}
Assume that $ F {\mathcal V }_{w_n} \not= 0 $ and 
consider the adjointness map 
$
\varphi = \varphi_{w_n} :  G \circ F  {\mathcal V }_{w_n}  \rightarrow {\mathcal V }_{w_n} $.
It is given concretely by multiplication
$$
\varphi: \, b_n e \otimes_{e b_n e} e  {\mathcal V }_{w_n}  \rightarrow {\mathcal V }_{w_n}, \,\, \, \, 
U \otimes e v \mapsto U e v.
$$
Using b) of the previous Lemma and the definition of left cells, we see that $ \varphi $ is surjective. 

\medskip
We now prove that $ \ker \varphi $ is zero. 
Recall from [MR] that 
$ U_i U_{i+1} \ldots U_{n-1 }  $, where $i = 0, 1, \ldots, n-1 $ generate $ b_n e $ 
as an $  e b_n e $-module. Using this and part a) of the previous Lemma 
we can write any $ k \in b_n e \otimes_{e b_n e} e  {\mathcal V }_{w_n}$ in the form  
$$ k = \sum_{i=0,1,\ldots, n-1} \, \sum_{  \,  w_{n-2} \in {\frak C}_{w_n} \cap W_{n-2} } 
\lambda_{i,w_{n-2}} \, U_i \, U_{i+1} \, \ldots \, U_{n-1} \otimes_{e b_n e} U_{n-1} C_{w_{n-2}} $$
where $ \lambda_{i,w_{n-2}} \in \C $.
Since $ U_{n-1} $ and $ C_{w_{n-2}} $ commute we have 
$$ U_{n-1} C_{w_{n-2}} = -\frac{1}{[2]} U_{n-1} C_{w_{n-2}} U_{n-1}  = 
-\frac{1}{[2]}\,  e \, C_{w_{n-2}} \, e . $$ 
Assume now that $ k \in \ker \varphi $. We then get 
$$
\begin{array}{c}
k= -\frac{1}{[2]}\,	
\sum_{i} \, \sum_{    w_{n-2}  } 
\lambda_{i,w_{n-2}} \, U_i \, U_{i+1} \, \ldots \, U_{n-1} \otimes_{e b_n e} 
e \, C_{w_{n-2}} \, e = 
\\
 -\frac{1}{[2]}\,  \sum_{i} \, \sum_{    w_{n-2}  } 
\lambda_{i,w_{n-2}} \, U_i \, U_{i+1} \, \ldots \, U_{n-1} U_{n-1} 
\, C_{w_{n-2}}  \otimes_{e b_n e} U_{n-1} 
\end{array}
$$
which is zero since $ k \in \ker \varphi$. This proves that $ \varphi $ is an isomorphism.

\medskip
Using a) of the previous Lemma once again, we now deduce that 
$$ F {\mathcal V }_{w_n} \simeq {\mathcal V }_{w_{n-2}} 
\, \,  \mbox{ for } \, w_n = {w_{n-2}}\,  s_{n-1}, \, \,  {w_{n-2}} \in  {\frak C}_{w_n} \cap W_{n-2}. $$
By Corollary 3.8 of [BGIL], $ {w_{n-2}} $ is independent of the choice of $ w_{n}$.
Under domino insertion, $ Sh(w_n) $ is obtained from $ Sh(w_{n-2}) $ by adding two horizontal dominoes, one
at the end of each line. Hence, using the formulas for $ \mathcal Q $ given in the proof of Theorem 
\ref{coideal}, we find that 
$$ {\mathcal Q}(Sh(w_{n-2})) = (a-1), (b-1 )  \, \, \, \mbox{if } \, \, {\mathcal Q}(Sh(w_{n})) = 
(a), (b ) $$ 
and hence the difference is the same. 

\medskip
If $ F {\mathcal V }_{w_n} = 0 $ we get by c) of the previous Lemma that 
$  {\mathcal V }_{w_n} \simeq \Delta_n(\pm n) $ and hence that $ \dim {\mathcal V }_{w_n} = 1 $. 
But then the combinatorial description of left cells in terms of domino tableaux gives 
$ w_n = 1 $ or $ w_n = s_0 $. For $ w_n =1 $ we have 
$ {\mathcal Q}(Sh(w_n ))= (n),(\emptyset)$ 
whereas for $ w_n = s_0 $ we have $ {\mathcal Q}(Sh(w_n ))= (\emptyset),(n)$, 
compatible with the actions of $ U_0 $ in ${\mathcal V }_{w_n} $.
The Theorem is proved.
\end{proof}
\begin{remark}
We think that the Theorem is valid for more general choices of $ q $ and $ Q$ within 
$ r = 0 $.
\end{remark}
\section{The Fock space}
In this section we give two applications of Theorem \ref{main} that both rely on 
the Fock space
approach to the representation theory of $ {\mathcal H}_n$. The first gives a new proof 
of the main results of [MW] using Ariki's Theorem and the second settles the question of 
determining the 
Kleshchev bipartition that corresponds to the simple $ b_n $-module 
$ L_n(\lambda) $. 
To set this up we first need the following 
Theorem.
\begin{theorem}{\label{Grothendieck}} 
In the Grothendieck group of $ b_n $-modules the equality $ \Delta_n( \lambda) 
=  S_n(a,b) $ holds where 
$  \lambda = a-b $
and $ S_n(a,b) $ is the Dipper-James-Murphy Specht module for $ {\mathcal H}_n $
corresponding to the bipartition $ (a),(b) $.
\end{theorem}
\begin{proof}
This follows basically from 
Theorem 3 and Theorem 6 of [RH]. 
On the other hand, since [RH] is based on a realization of $ b_n $ as a quotient of the Ariki-Koike 
algebra $ {\mathcal AK}_n(\lambda_1, \lambda_2,q )  $ 
and a realization of $ S_n(a,b) $ as a permutation module in the Ariki-Yamada-Terasoma tensor space 
for $ {\mathcal AK}_n(\lambda_1, \lambda_2,q )  $, 
we still give a few details on how to convert from one situation to the other.



\medskip
Recall that $ TLB_{n} = {\mathcal H}_{n} / J_{n} $ is free over $ \mathcal A $ and so we have that
$$
TLB_{n, q, Q} = {\mathcal H}_{n, q, Q} / J_{n,q,Q}.
$$
By Lemma \ref{relates} we also know that $ TLB_{n, q, Q} $ 
is isomorphic to the blob algebra $ b_n( q, m) $, that on the other hand 
was realized in [RH] as $ {\mathcal AK}_n(\lambda_1, \lambda_2,q ) / G_n $
where $ \lambda_1 = \frac{q^m}{q-q^{-1}}, \lambda_2 = \frac{q^{-m}}{q-q^{-1}} $
and $ G_n  $ is the ideal of $ {\mathcal AK}_n(\lambda_1, \lambda_2,q )  $ 
generated by $ (X_1 X_2 - \lambda_1 \lambda_2)(g_1 -q ) $. 
The last realization requires the conditions 
$ q^4 \neq  1, \lambda_1 \neq  \lambda_2, \lambda_1 \neq q^ 2 \lambda_2 $ and these
conditions are imposed throughout [RH]. 

\medskip
Instead of 
converting directly between the Hecke algebra and the Ariki-Koike algebra setting
we prefer to proceed as follows.

\medskip
The Hecke algebra ${\mathcal H}_n $ {\it is} an Ariki-Koike algebra
with parameters $ \lambda_1 = Q, \lambda_2 = -1/Q $ and so we can develop the theory of [RH] entirely from 
the ${\mathcal H}_n $ point of view, once we have
proved that $ J_n $ acts trivially in the Ariki-Terasoma-Yamada tensor space $ V^{ \otimes n} $ 
for these choices of parameters
when $ \dim V = 2 $, corresponding to Theorem 1 of [RH].
 
\medskip 
Let us therefore detail how the analogue of Theorem 1 of [RH] is proved.
Let $ V$ be a complex vector space of dimension two and let $ v_1, v_2 $ be a basis. 
Let 
$ R \in \End_{\C } (V \otimes V) $ be given by 
$$   \begin{array}{l} R( v_i \otimes v_j ) =   
q v_i \otimes v_j  \,\,\,\,\,\,\,\,\,\,   \mbox{if } i=j  \\
R( v_2 \otimes v_1 ) = 
v_1 \otimes v_2 \\ R( v_1 \otimes v_2 ) = 
v_2 \otimes v_1 + (q-q^{-1})  v_1 \otimes v_2. 
\end{array} 
$$
For $ i = 1, 2, \ldots ,n-1 $, we let $ T_i \in {\mathcal H}_n $ act in the tensor space 
$ V^{\otimes n } $ by
$$  T_{i} :=Id^{\otimes i-1 } \otimes R \otimes Id^{\otimes n-i-1} .$$
For 
$ v = v_{i_1} \otimes v_{i_2} \otimes \ldots \otimes v_{i_k}   \otimes v_{i_{k+1}} \otimes \ldots \otimes v_{i_n } $,
we define $ S_k \in \End_{\C } (V^{\otimes n} ) $ by 
$$ S_k( v ) = \left\{
\begin{array}{ll} q v_{i_1} \otimes v_{i_2} \otimes \ldots \otimes v_{i_{k}} \otimes v_{i_{k+1}}    \otimes \ldots \otimes v_{i_n}  
& \mbox{if   }   i_k = i_{k+1} \\
v_{i_1} \otimes v_{i_2} \otimes \ldots \otimes v_{i_{k+1}} \otimes v_{i_{k}} \otimes \ldots \otimes v_{i_n}   & \mbox{otherwise }
\end{array} \right\}
$$
and let $ \varpi \in \End_{\C } (V^{\otimes n} ) $ be given 
by
$$ \varpi(v_{i_1} \otimes  
 v_{i_2} \otimes   \ldots \otimes v_{i_n } ) := \left\{ 
\begin{array}{ll} Q v & \mbox{ if } i_1 = 1 \\
-Q^{-1} v & \mbox{ if } i_1 = 2.
\end{array} \right. 
$$
Setting  
$$ T_0 := T_{1}^{-1} \ldots T_{n-2}^{-1} T_{n-1}^{-1} S_{n-1} S_{n-2} \ldots S_1 \varpi $$   
it is then proved in [ATY] that the linear maps 
$ T_0, T_1, \ldots, T_{n-1} $ define an action of $ {\mathcal H}_n $ in $ V^{\otimes n}$.
Endowed with this action of $ {\mathcal H}_n $, we call $ V^{\otimes n}$ the Ariki-Yamada-Terasoma tensor space.

\medskip
Let us now show that the ideal $ J_n $ is annihilated under this action. This 
is well-known for the generator 
$ C_1 C_2 C_1 - C_1 $ so we concentrate on
$ C_1 C_0 C_1 - [2]_{\frac{Q}{q}} C_1 $.  
Since $ C_1 $ acts semisimply in $ \Span \{ \,  v_i \otimes v_j \, | \, i, j= 1,2 \, \} $ with 
eigenvalue $ 0 $ of multiplicity three and eigenvalue $ -[2] $ of multiplicity one, it 
is enough to check the relation on vectors of the form $ C_1 v $ where 
$ v = v_{2} \otimes v_{1} \otimes v_{i_3} \otimes    \ldots \otimes v_{i_n } $
since $ C_1 v \not= 0 $ for such $ v$. But $ C_1 v = (v_1 \otimes v_2 -q v_2 \otimes v_1) \otimes
v_{i_3} \otimes    \ldots \otimes v_{i_n } $ is an eigenvector for $ C_1 $ of eigenvalue 
$ -[2] $ and hence it is enough to show that 
\begin{equation}{\label{ideal_vanish}}
 C_1 C_0  (v_1 \otimes v_2 -q v_2 \otimes v_1) \otimes
\overline{v}
= [2]_{\frac{Q}{q}} (v_1 \otimes v_2 -q v_2 \otimes v_1) \otimes
\overline{v}  
\end{equation}
where $\overline{v} = v_{i_3} \otimes    \ldots \otimes v_{i_n } $.
Let us consider the left hand side of this equation. Using Lemma 1 of [RH], which is a reformulation of a result 
of [ATY],
we find that $$ C_1 C_0 \, q \, v_2 \otimes v_1 \otimes \overline{v} = 
q^2 (  Q + Q^{-1}) \,C_1 v_1 \otimes v_2 \otimes \overline{v} .$$
We then consider 
$ C_1 C_0 \, v_1 \otimes v_2  \otimes
\overline{v} $
which we rewrite as follows 
$$ 
\begin{array}{r}
C_1 C_0 \, v_1 \otimes v_2  \otimes
\overline{v} =  (T_1 - q)  (T_0 -Q)  \, v_1 \otimes v_2  \otimes
\overline{v} =  \\
(T_1 - q)  T_0   \, v_1 \otimes v_2  \otimes
\overline{v} - Q  
C_1   \, v_1 \otimes v_2  \otimes
\overline{v}.
\end{array}
$$
We here consider the first term $ (T_1 - q)  T_0   \, v_1 \otimes v_2  \otimes
\overline{v} $ which we rewrite as follows 
$$ 
\begin{array}{r}
(T_1 - q)  T_0   \, v_1 \otimes v_2  \otimes
\overline{v}
=
-q (T_1 - q) T_1   T_0 T_1  \, v_2 \otimes v_1  \otimes
\overline{v} = \\
-qQ (T_1 - q)   \, v_2 \otimes v_1  \otimes
\overline{v} 
= q^2 Q C_1  \, v_1 \otimes v_2  \otimes
\overline{v}
\end{array}
$$
where we for the second equality used the argument given in the proof of Theorem 1 of [RH].
Summing up, the LHS of (\ref{ideal_vanish}) equals
$$ (-Q -q^2 Q^{-1}) C_1 v_1 \otimes v_2  \otimes
\overline{v} $$
which coincides with the RHS.

\medskip
We can now develop the theory of [RH] from the Hecke algebra point of view. 
Especially, for $ \lambda \in \Lambda_n $ we define the permutation module
$$ M_n(\lambda) :=  \Span_{\C} \{ v_{i_1} \otimes v_{i_2} \otimes \ldots \otimes v_{i_n} \, | \,
\,  \#\{k :  i_k = 1\} -  \#\{k :  i_k = 2\}   = \lambda \,\} $$
and get that $ M_n(\lambda) $ satisfies the functorial properties for $ F $ of 
(\ref{functors}).

\medskip
Theorem 3 of [RH] is proved by induction. One checks that the inductive step works for all choices
of the parameters satisfying $ \lambda_1 \neq \lambda_2 $. 
But $ \lambda_1 = Q = i q^m  $ and $ \lambda_2 = -Q^{-1} = i q^{-m} $ and so we have 
$ \lambda_1/ \lambda_2 = q^{2m} = -q  \neq 1$, as needed. The induction basis is based on Lemma 3 of [RH].
The proof of that Lemma works provided that 
$ \lambda_1 (q-q^{-1} ) \neq q( \lambda_1 - \lambda_2) $. But this is equivalent to
$ -q \neq q^2 $ that is $ q \neq -1 $, as needed.


\medskip
Finally the proof of Theorem 6 of [RH] 
claiming that 
$ M_n (\lambda) \cong S_n( a, b )^{\circledast}  $ 
is independent of the choices of the parameters and 
goes directly over. But in the Grothendieck group of $ b_n$-modules,  
$ S_n( a, b ) $ is equal to its contragredient dual $S_n( a, b )^{\circledast} $, 
and so the proof of the Theorem is finished.
\end{proof}
\begin{remark}
In view of Theorem \ref{main}, an alternative proof might have been 
obtained using the results of section 4 of [P1]. 
\end{remark}

\begin{remark}
At this point we may remark that combining  
Theorem {\ref{main}} with Lemma 2 of [RH], we get many examples of cells 
modules for different choices of $ r $ that are not isomorphic. 
Indeed Lemma 2 of [RH] gives many examples of the adjointness map 
$ G \circ F M_n( \lambda  ) \rightarrow M_n( \lambda ) $ failing to be an isomorphism.
Note
that the condition in that Lemma 2, that $ q$ be an odd order root of unity, is not 
needed for showing that the adjointness map is not surjective -- as is indeed mentioned in the proof of that Lemma 2.
\end{remark}
We now recall the Fock space approach to the representation theory of $ {\mathcal H}_n$. 
Let $ s= (s_1, s_2 ) \in {\Z}^2 $ and let $ {\frak F}^s $ be the associated Fock space 
of level two. As a $ {\C}(v)$-vector space it is given by 
$$ {\frak F}^s = \bigoplus_{\lambda \in \Bip} {\C}(v) \, | \,\lambda, s \rangle $$
where $ | \lambda, s \rangle $ is a symbol. 
Set  $ e := l/2 $ where $ l $ is the multiplicative order of $ q$ as in ({\ref{condition}}).
Let us briefly recall how $  {\frak F}^s $  
becomes an
integrable module for 
the quantum group $ {\mathcal U}_v ( \widehat{\frak sl}_e) $, following the construction in [JMMO]. 

\medskip
Since $ {\mathcal U}_v ( \widehat{\frak sl}_e) $ is the $ \C(v) $-algebra generated by 
$ e_i, f_i, \, i=0,1, \ldots, e-1 $ and $ k_h, \, h \in \frak h $ subject to certain well-known 
relations, it is enough to explain how these generators act in 
$ {\frak F}^s $. 
To any bipartition $ (\lambda^{(1)}, \lambda^{(2)}) $ we associate its diagram 
$$ \{ (i,j,c) \, |\,  c=1,2 \mbox{ and }1 \leq j \leq \lambda_i^{(c)} \,\}. $$
For a node $ \gamma =(i,j,c) $ of $ (\lambda^{(1)}, \lambda^{(2)}) $ we define its $e$-residue 
by $ \Res_e(\gamma)= j-i + s_c \mod e $. 
We define a total order on the nodes of $ (\lambda^{(1)}, \lambda^{(2)}) $ 
by $ \gamma = (i,j,c) < \gamma^{\prime} = (i^{\prime},j^{\prime},c^{\prime})$ if 
$ j-i + s_c < j^{\prime}- i^{\prime}+ s_c^{\prime} $ or if 
$ j-i + s_c = j^{\prime}- i^{\prime}+ s_c^{\prime} $ and 
$  c^{\prime} < c $ (notice this last inequality!).
If $ \lambda = (\lambda^{(1)}, \lambda^{(2)}) $ and 
$ \mu = (\mu^{(1)}, \mu^{(2)}) $ are bipartitions such that $ \lambda \subset \mu $ and 
$ \gamma = \mu \setminus \lambda $ is an $ i $-node we say that $ \gamma $ is a removable 
$ i$-node of $ \mu $ and an addable $ i$-node of $ \lambda $ and we set 
$$ \begin{array}{rl}
N_i^{>}(\lambda, \mu) := & \# \{ \, \mbox{ addable } i\mbox{-nodes} \, \gamma^{\prime}  \mbox{ of } 
\lambda \mbox{ such that } \gamma^{\prime} > \gamma \, \} \\ 
- & 
\# \{ \, \mbox{ removable } i\mbox{-nodes} \,  \gamma^{\prime}  \mbox{ of } 
\lambda \mbox{ such that } \gamma^{\prime} > \gamma \, \} \\
N_i^{<}(\lambda, \mu) := & \# \{ \, \mbox{ addable } i\mbox{-nodes} \, \gamma^{\prime}  \mbox{ of } 
\lambda \mbox{ such that } \gamma^{\prime} < \gamma \, \} \\ 
- & 
\# \{ \, \mbox{ removable } i\mbox{-nodes} \,  \gamma^{\prime}  \mbox{ of } 
\lambda \mbox{ such that } \gamma^{\prime} < \gamma \, \}. \\
\end{array} 
$$
The actions of $ f_i, e_i $ on a basis vector of $ {\frak F}^s $ are now as follows
$$ \begin{array}{l} 
 f_i \, |\lambda, s \rangle = \sum_{ \mu, \, \Res(\mu \setminus \lambda) \equiv i \, \Mod  e } 
\, v^{ N_i^{>}(\lambda, \mu)} \, | \mu,s \rangle \\
 e_i \, |\mu, s \rangle = \sum_{ \lambda, \, \Res(\mu \setminus \lambda) \equiv i \,  \Mod  e} 
\, v^{- N_i^{<}(\lambda, \mu)} \, | \lambda,s \rangle. 
\end{array}
$$
There are similar formulas for the other generators.
It is one of the important issues of the Fock space approach to the representation theory 
of $ {\mathcal H}_n $ that 
$ {\frak F}^s $ with this action 
not only depends on the 
classes $ s_1 \, \mbox{mod} \,e  $ and $  s_2 \, \mbox{mod} \, e $, 
but on $ s $ itself.

\medskip
Let $ {\mathcal U}_v ( \widehat{\frak sl}_e) \rightarrow {\mathcal U}_v ( \widehat{\frak sl}_e), \,
u \mapsto \overline{u} $ be the bar involution given by 
$$ \begin{array}{cccc}
\overline{v}:= v^{-1}, & 
\overline{f_i}:= f_i, &
\overline{e_i}:= e_i, & 
\overline{k_h}:= k_{-h} 
\end{array}
$$
and let $ {\frak F}^s \rightarrow {\frak F}^s, \,  x \mapsto \overline{x} $
be the bar involution of the Fock space constructed by Uglov in [U].
It satisfies $ \overline{ \emptyset, s \rangle } = \emptyset, s \rangle $ and
is compatible with the bar involution on 
$ {\mathcal U}_v ( \widehat{\frak sl}_e) $, that is 
$ \overline{ u x } = \overline{u} \,  \overline{x} $ for $ u \in {\mathcal U}_v ( \widehat{\frak sl}_e) $
and $ x \in {\frak F}^s $. By the results of [U] we get 
for $ \lambda \in \Bip $ a unique $ G(\lambda, s) \in {\frak F}^s $ such 
that 
$$ 
\begin{array}{cc}
\overline{G(\lambda, s)} = G(\lambda, s),  & 
G(\lambda, s) \equiv | \lambda, s \rangle \mod v \C[v] \,{\frak F}^s.
\end{array}
$$
Write for $ \mu \in \Bip $ 
$$ G(\mu, s) = \sum_{\lambda \in \Bip } d_{\lambda, \mu }^s (v) | \lambda, s \rangle. $$
Set $ {\mathcal M}[s] :=  {\mathcal U}_v ( \widehat{\frak sl}_e) \, | \, \emptyset, s \rangle $. 
Then $ {\mathcal M}[s] $ is an integrable module for $ {\mathcal U}_v ( \widehat{\frak sl}_e) $ and
so the crystal/canonical basis theory applies to it. In fact, 
there is a subset $ \Bip_e^s \subset \Bip $ such that $ G(\lambda, s) $ for $ (\lambda, s) \in \Bip_e^s $ is the 
canonical basis/global crystal basis of ${\mathcal M}[s] $. Set $ \Bip_e^s(n) := \Bip_e^s \cap \Bip(n)$. 
Assume that $ m \equiv s_1 -s_2 $. Then by the deep Theorem of Ariki in [A], we have that   
$ \Bip_e^s(n) $ parameterizes the irreducible modules for $ {\mathcal H}_n$ 
with corresponding decomposition numbers $ d_{\lambda, \mu }^s (1) $.

\medskip

The proof of our next Theorem is 
essentially the same as the proof of Theorem 4.7 of [BJ], but notice that 
Theorem 4.7 of [BJ] requires the validity of the Conjectures A, B and $\rm{B}^{\prime}$ of [BJ].
As already mentioned,  Conjecture A holds in the $ r = 0 $ case whereas, as we shall see, 
we can replace Conjecture B by Theorem \ref{main} and Conjecture $\rm{B}^{\prime}$ by our previous
Theorem \ref{Grothendieck}.

\begin{theorem}{\label{Bonnafe-Jacon}} 
Let $ m,l, e $ be as above and let $ p $ be the largest integer such that $ m + pe \leq 0 $ and 
set $ s:=(m+pe,0) $. Then for $ \mu \in \Bip_e^s(n) $ we have 
\begin{equation}{\label{crystal_order}} G(\mu, s) = | \mu, s \rangle + \sum_{\lambda \in \Bip(n),\, \lambda \prec \mu } d_{\lambda, \mu }^s (v) | \lambda, s \rangle.
\end{equation}
Moreover, identifying $ \tau =(t_1, t_2) \in \Bip_1(n) $ with 
$ f(\tau)=t_1-t_2 \in \Lambda_n$ 
we have for $ \lambda, \mu \in \Bip_1 $ that 
\begin{equation}{\label{FLOTW}}
[\Delta_n(\lambda), L_n(\mu)]
= d_{\lambda, \mu}^s(1).
\end{equation}
\end{theorem}
\begin{proof}
By the choice of $ s $ we have 
formula (\ref{crystal_order}) as in the proof of Theorem 4.7 of [BJ]. 
Notice now that $ m + pe \not=0 $. 
Thus we have that 
$ \lambda $ and $ \mu $ of 
(\ref{FLOTW}) belong to $ \Bip_e^s(n) $. These bipartitions were studied in [FLOTW] 
and are called FLOTW bipartitions in [Ja, BJ] and other references.

\medskip
Take now $ \nu = (n_1, n_2)  \in \Bip_1(n) $ corresponding to $ \nu \in \Lambda_n $.
According to Ariki's Theorem there exists $ \mu \in \Bip_e^s(n) $ such that the 
decomposition number $ d_{\lambda, \nu} := [ S_n(\lambda), L_n(\nu)] $ satisfies
$$[ \Delta_n(\lambda), L_n(\nu)] =  d_{\lambda, \nu}  = d_{\lambda, \mu}^s(1) $$
for all $ \lambda \in \Bip(n) $ where we used Theorem 
\ref{Grothendieck}
for the first equality.
Setting $ \lambda = \nu $ we get that $ \nu \preceq \mu $ and 
setting $ \lambda = \mu $ we get that $ \mu \preceq \nu $. Hence $ \mu = \nu $ and the 
Theorem is proved.
\end{proof}

The next step is now to calculate the numbers $ d_{\lambda, \mu}^s(1) $ for
$ \lambda, \mu \in \Bip(n) $. Uglov's proof of the existence of $ G(\lambda,s) $ 
is not straightforward, but still
constructive.  Notice that his algorithm has been simplified by 
Yvonne in [Y], and that Fayers, [Fa], as well as Jacon, [Ja1], have found combinatorial 
generalizations of the LLT algorithm to higher levels.

\medskip
On the other hand, since we here only focus on bipartitions in $ \Bip_1(n)$ actually the 
properties of $ G(\lambda,s) $ already mentioned are sufficient to calculate $ G(\lambda,s) $ and hence
$ d_{\lambda, \mu}^s(1) $. 

\medskip
Indeed, 
set $ m_- :=-(m+(p+1)e) $ and 
recall from [MW] that the choices of $ e $ and $ m $ determine an alcove geometry in $ \R $
with zero dimensional walls in the integral points $ {\mathcal M} :=\{  m_- + k e \, | \, k \in \Z \}$ and 
fundamental alcove $ A_0 $ being the one that contains $ 0 $. The associated Weyl group $ \mathcal W $ is 
infinite dihedral, generated by $ s_+ $ and $ s_-$ where $ s_+ $ ($s_-$) is the reflection in 
the right (left) wall of the fundamental alcove. Set 
$ \Lambda_n^{reg}:= \Lambda_n \setminus {\mathcal M} $ and for 
$ \lambda \in \Lambda_n^{reg} $ write 
$ A_{\lambda} $ for 
the alcove containing $ \lambda $.
For $ \lambda \in \Lambda_n^{reg} $ we 
define $ w_{\lambda} \in \mathcal W $ by the condition 
$ w_{\lambda } A_0 = A_{\lambda} $. Thus $ w_{\lambda } < w_{\mu  } $ in the 
Bruhat-Chevalley order implies $ \lambda > \mu $ in the quasi-hereditary order.
We can now formulate the next Theorem. 
The second part of it was proved in [MW] using completely different methods.
\begin{theorem}{\label{Decomposition numbers}} 
Let $ \lambda, \mu  \in \Lambda_n^{reg}$. Then we have 
\begin{equation}{\label{decom1}} 
d_{\lambda, \mu}^s(v) = \left\{ \begin{array}{cc} v^{ l(w_{\lambda}) -  l(w_{\mu}) } & \mbox{ if }  
w_{\lambda} \leq w_{\mu} \\ 0 & \mbox{ otherwise } \end{array} \right.
\end{equation}  
\begin{equation}{\label{decom2}}  
[ \Delta_n(\lambda), L_n(\mu)] = 
\left\{ \begin{array}{cc} 1 & \mbox{ if }  
w_{\lambda} \leq w_{\mu} \\ 0 & \mbox{ otherwise. } \end{array} \right.
\end{equation}  
\end{theorem}
\begin{proof}
Following [MW] we enumerate the elements of $ \mathcal W $ as follows 
$$
w_i = \left\{ \begin{array}{lr} 1 & \mbox{ if } i = 0 \\
s_- s_+ s_i \ldots ( \, -i \mbox{ terms} ) & \mbox{ if } i < 0 \\
s_+ s_- s_+ \ldots ( \, i \mbox{ terms} ) & \mbox{ if } i > 0 
\end{array} \right.
$$
and define $ A_i := w_i A_0 $. Then, $ A_i $ is the alcove at distance $ i $ from $ A_0 $, 
positioned to the right if $ i $ is positive and to the left if $ i$ is negative.

\medskip
Write $ s_1 := m+pe $ such that $ s= (s_1, 0 ) $. Set furthermore
$ m_+ := m_- +e $. Then the fundamental alcove is limited by $ m_- $ and $ m_+$.
Assume now that $ \lambda= (k_1,k_2) $ belongs to $ A_i \cap \Lambda_n^{reg}$ with $i \geq 0$.
Let 
$ r_1, r_2 \in \{0, 1 , \ldots , e-1 \} $ be the residues modulo $ e $ of $ k_1+s_1, k_2 $. 

\medskip
We now act with elements of the form 
$  f_{r_1+j}  \ldots  f_{r_1+1}  f_{r_1} $ in $ | \lambda, s \rangle $ 
and consider the images in $ {\frak F}^{s, \, \geq 2} :=    {\frak F}^s / I^{ \geq 2} $ 
where $ I^{ \geq 2} := \Span\{ | \nu, s \rangle \, | \, \nu \notin \Bip_1 \} $.
These images 
move towards the right wall of $ A_i $. The wall will be reached when $ r_1 + j = r_2 \, \mbox{ mod } \, e $
and the image will be $  | \mu , s \rangle $ where $ \mu = (k_1 + r_1 - r_2, k_2 )  $, 
i.e. with 
$ v $ power equal to $ v^0 $ since $  k \geq 0 $. Notice here that 
the wall $ m_+ $ of $ A_0 $ corresponds exactly to the second case in the definition of 
the order relation on the nodes. 

\medskip 
In the formalism of translation functors, as 
exposed for example in [S], the 
process just described corresponds to translation upwards on the wall. 

\medskip 
Acting with $ f_{r_1 } $ in $  | \mu , s \rangle $ and considering the images in 
$ {\frak F}^{s, \, \geq 2} :=    {\frak F}^s / I^{ \geq 2} $ the result is 
$$ | \mu^{up}, s \rangle + v \, | \mu^{down }, s \rangle $$
where $ \mu^{up} =   (k_1+1, k_2 + r_1 - r_2 ) $ and $ \mu^{down} =   (k_1, k_2 + r_1 - r_2 +1) $
and once again we 
get correspondence with the translation functor formalism. 

\medskip
Similarly, we go through the other cases and find that translation upwards through the wall behaves 
as above whereas translation downwards through the wall 
$ | \mu, s \rangle  $ 
is given by  
$$ v^{-1} | \mu^{up}, s \rangle + \, | \mu^{down }, s \rangle $$
where $ \mu^{up} $ and $ \mu^{down}$ are chosen analogously to the first case.

\medskip
Using these rules, together with (\ref{crystal_order}) and Theorem \ref{coideal} it us 
now straightforward to 
calculate $ G(\lambda, s) $ modulo $ I^{ \geq 2} $ for $ \lambda \in \Bip_1 $ to obtain
formula (\ref{decom1}).
Finally, formula (\ref{decom2}) then follows from the previous Theorem.
\end{proof}

\medskip
Let us finish by mentioning one more application of our results. 
The (negative) Kleshchev bipartitions, see [Ja], [BJ] and references therein, are by definition 
the bipartitions that belong to $ \Bip_e^s $ where $ s= (d+re, 0 ) $ is chosen such that 
$ d+re > n-1-e $, this is the so-called (negative) asymptotic case. 
The Kleshchev bipartitions parameterize
the simple modules for ${\mathcal H}_n $
when we use the dual Specht modules to parameterize, that is when we use the $M_n(\lambda) $'s.

\medskip
The question raised in [RH] of
determining the Kleshchev bipartition $ \lambda =(l_1, l_2) $ that corresponds 
to the simple $ b_n $-module with parameter $ \tau=(t_1, t_2) $ can now be 
solved by applying Kashiwara's operators to the crystal graphs of the Fock spaces. 
Consider as an example $ e = 3, m = 2 $. Then $ s= (-1, 0 ) $. 
In the crystal graph of $ {\mathcal M } ( -1, 0) $ we have
$$ \tilde{f}_0 \tilde{f}_1 \tilde{f}_0 \tilde{f}_2\tilde{f}_2 \tilde{f}_1 \tilde{f}_1 \tilde{f}_0 \tilde{f}_0 
\tilde{f}_2  (\emptyset, \emptyset ) = (6, 4) $$
whereas the same sequence of crystal operators sends $ (\emptyset, \emptyset ) $ to $ ((6,3), (1)) $ in 
$ {\mathcal M } ( 11, 0) $.

\medskip
Jacon has constructed in [Ja] an algorithm for converting between such crystal graphs. The following tables
have been calculated using an 
implementation of his algorithm in the GAP system. They convert between the bipartitions in 
$ \Bip_1(10) $ and the corresponding Kleshchev bipartitions, that we denote $ \KBip_1(10) $.

$$
\begin{array}{lccr}
\begin{tabular}[t]{|c|c|} 
\hline
\multicolumn{2}{|c|}{$ e = 3, m = 2, s= (-1, 0 )$} \\
\hline 
$ \Bip_1(10)$ & $ \KBip_1(10)$\\
\hline
$ (10),(\emptyset)$ & $(10), (\emptyset)$ \\
\hline
$ (9),(1)$ & $(9),(1)$ \\
\hline
$ (8),(2)$ & $((8,1), (1))$ \\
\hline
$ (7),(3)$ & $((7,2), (1))$ \\
\hline
$ (6),(4)$ & $((6,3), (1))$ \\
\hline
$ (5),(5)$ & $((5,4), (1))$ \\
\hline
$ (4),(6)$ & $((4,2), (4))$ \\
\hline
$ (3),(7)$ & $((5,1), (4))$ \\
\hline
$ (2),(8)$ & $(6),(4)$ \\
\hline
$ (1),(9)$ & $(7),(3)$ \\
\hline
$ (\emptyset),(10)$ & $(8), (2)$ \\
\hline
\end{tabular}
&  & & 
\begin{tabular}[t]{|c|c|} 
\hline
\multicolumn{2}{|c|}{$ e = 5, m = 3, s= (-2, 0 )$} \\
\hline
$ \Bip_1(10)$ & $ \KBip_1(10)$\\
\hline
$ (10),(\emptyset)$ & $(10), (\emptyset)$ \\
\hline
$ (9),(1)$ & $(9),(1)$ \\
\hline
$ (8),(2)$ & $(8), (2)$ \\
\hline
$ (7),(3)$ & $((7,1), (2))$ \\
\hline
$ (6),(4)$ & $((6,2), (2))$ \\
\hline
$ (5),(5)$ & $((5,3), (2))$ \\
\hline
$ (4),(6)$ & $((4,4), (2))$ \\
\hline
$ (3),(7)$ & $(4), (6)$ \\
\hline
$ (2),(8)$ & $(5),(5)$ \\
\hline
$ (1),(9)$ & $(6),(4)$ \\
\hline
$ (\emptyset),(10)$ & $(7), (3)$ \\
\hline
\end{tabular}
\end{array}
$$

$$
\begin{array}{lccr}
\begin{tabular}[t]{|c|c|} 
\hline
\multicolumn{2}{|c|}{$ e = 7, m = 4, s= (-3, 0 )$} \\
\hline 
$ \Bip_1(10)$ & $ \KBip_1(10)$\\
\hline
$ (10),(\emptyset)$ & $(10), (\emptyset)$ \\
\hline
$ (9),(1)$ & $(9),(1)$ \\
\hline
$ (8),(2)$ & $(8), (2)$ \\
\hline
$ (7),(3)$ & $(7), (3)$ \\
\hline
$ (6),(4)$ & $((6,1), (3))$ \\
\hline
$ (5),(5)$ & $((5,2), (3))$ \\
\hline
$ (4),(6)$ & $((4,3), (3))$ \\
\hline
$ (3),(7)$ & $(3), (7)$ \\
\hline
$ (2),(8)$ & $(4),(6)$ \\
\hline
$ (1),(9)$ & $(5),(5)$ \\
\hline
$ (\emptyset),(10)$ & $(6), (4)$ \\
\hline
\end{tabular}
&  & & 
\begin{tabular}[t]{|c|c|} 
\hline
\multicolumn{2}{|c|}{$ e = 9, m = 5, s= (-4, 0 )$} \\
\hline
$ \Bip_1(10)$ & $ \KBip_1(10)$\\
\hline
$ (10),(\emptyset)$ & $(10), (\emptyset)$ \\
\hline
$ (9),(1)$ & $(9),(1)$ \\
\hline
$ (8),(2)$ & $(8), (2)$ \\
\hline
$ (7),(3)$ & $(7), (3)$ \\
\hline
$ (6),(4)$ & $(6), (4)$ \\
\hline
$ (5),(5)$ & $((5,1), (4))$ \\
\hline
$ (4),(6)$ & $((4,2), (4))$ \\
\hline
$ (3),(7)$ & $((3,3), (4))$ \\
\hline
$ (2),(8)$ & $(3),(7)$ \\
\hline
$ (1),(9)$ & $(4),(6)$ \\
\hline
$ (\emptyset),(10)$ & $(5), (5	)$ \\
\hline
\end{tabular}
\end{array}
$$

\medskip \medskip 

It can be seen that 
the correspondence between $ \Bip_1(10)$ and $ \KBip_1(10)$ 
works as the identity 
in the top $ m$ lines of all of these tables. 
This is a consequence of the following last result. 

\begin{Prop}{\label{first_m}} 
For $ n_2 = 0, 1, \ldots ,m-1 $ we have that $ \Delta_n(\lambda) \cong  S_n(n_1, n_2 )^{\circledast}  $ 
where $ \lambda:=n_1-n_2$.
In particular, for these values of $ n_2 $ 
the irreducible $ b_n $-module $ L_n(\lambda)    $ is isomorphic to the unique irreducible quotient of 
$ S_n(n_1, n_2 )^{\circledast} $.
\end{Prop}
\begin{proof}
We already saw that $ S_n(n_1, n_2 )^{\circledast}  \cong M_n(n_1, n_2) $.
Let $ N $ be the minimum of $ n_1 $ and $ n_2$. Then we get by the 
categorical properties ({\ref{functors}}), that
$$  \Delta_n(\lambda) \cong G^{ \circ N} \circ  F^{ \circ N} \Delta_n(\lambda) $$ 
where $ F^{ \circ N} :=F \circ \ldots \circ F $ and $ G^{ \circ N} :=G \circ \ldots \circ G $ are the $N$-fold compositions 
of the functors $ F $ and $ G$.
On the other hand, using Lemma 2 of [RH], we find for the chosen values of $ n_2 $, a similar description of 
$ M_n(n_1, n_2)$ as follows
$$  M_n(n_1, n_2)  \cong G^{ \circ N} \circ  F^{ \circ N}  M_n(n_1, n_2) . $$ 
Notice that the conditions of that Lemma on $ l $ to be odd and $ n_2 \neq m $ mod $ l $ can be replaced by 
$ n_2 \neq m $ mod $ e $, as can easily be seen from the proof of the Lemma.
We now use Theorem 2 and Lemma 3 of {\it loc. cit.} together with ({\ref{functors}}) 
to deduce that 
$$ F^{ \circ N} \Delta_n(\lambda)  \cong F^{ \circ N} M_n(n_1, n_2) $$ 
and the Proposition follows by applying $G^{ \circ N} $ to this isomorphism.
\end{proof}

\medskip
Unfortunately, in general we do not have a non-recursive description of the elements of $ \KBip_1(n) $.

\medskip \medskip

\bibliographystyle{amsplain}

\end{document}